\documentclass[11pt]{article}

\usepackage{amsmath, amsfonts, amssymb, amsthm,  graphicx, enumerate}
\usepackage[letterpaper, left=1truein, right=1truein, top = 1truein, bottom = 1truein]{geometry}
\usepackage[numbers, square]{natbib}

\usepackage[%
            CJKbookmarks=true,
            bookmarksnumbered=true,
            bookmarksopen=true,
            colorlinks=true,
            citecolor=red,
            linkcolor=blue,
            anchorcolor=red,
            urlcolor=blue,
            ]{hyperref}

\usepackage[usenames,dvipsnames]{color}

\usepackage[ruled, vlined, lined, commentsnumbered]{algorithm2e}

\usepackage{prettyref,soul,xspace}

\usepackage{tikz}
\usepackage{pgfplots,makecell}

\usepackage{bm}
\usepackage{fancyhdr}

\theoremstyle{remark}

\theoremstyle{plain}

\newtheorem{lemma}{Lemma}
\newtheorem{theorem}{Theorem}
\newtheorem{proposition}{Proposition}

\newcommand{\p}{\mathbb{P}}
\newcommand{\E}{\mathbb{E}}

\newcommand{\iid}{\stackrel{iid}{\sim}}
\newcommand{\dequal}{\stackrel{d}{=}}
\newcommand{\br}[1]{\left( #1 \right)}

\newcommand{\cbr}[1]{\left\{ #1 \right\}}
\newcommand{\pbr}[1]{\p\left( #1 \right)}
\newcommand{\ebr}[1]{\exp\left( #1 \right)}
\newcommand{\abs}[1]{\left| #1 \right|}

\newcommand{\mathr}{\mathbb{R}}
\newcommand{\mathn}{\mathcal{N}}
\newcommand{\mathf}{\mathcal{F}}

\newcommand{\indic}[1]{{\mathbb{I}\left\{{#1}\right\}}}
\newcommand{\iprod}[2]{\left \langle #1, #2 \right\rangle}
\newcommand{\norm}[1]{\left\|{#1} \right\|}
\newcommand{\normt}[1]{\|{#1} \|}
\newcommand{\opnorm}[1]{\norm{#1}_{\rm op}}
\newcommand{\opnormt}[1]{\normt{#1}_{\rm op}}
\newcommand{\opnorminf}[1]{\norm{#1}_{\rm op,\infty}}

\newcommand{\fnorm}[1]{\norm{#1}_{\rm F}}
\newcommand{\argmin}{\mathop{\rm argmin}}
\newcommand{\argmax}{\mathop{\rm argmax}}

\newcommand{\define}{:=}
\newcommand{\definei}{=:}

\newcommand{\hatH}{M}

\usepackage{subcaption}
\usepackage{bbm}

\newcommand{\loo}[1]{{{#1}^{(j)}}}
\newcommand{\calE}{\mathcal{E}}
\newcommand{\MN}{\mathcal{MN}}

\newcommand{\calN}{\mathcal{N}}
\renewcommand{\vec}{\text{vec}}
\newcommand{\zero}{0}

\newcommand{\R}{\mathbb{R}}
\usepackage{import}
\newcommand{\inner}[2]{\left\langle {#1}, {#2}\right\rangle}
\newcommand\numberthis{\addtocounter{equation}{1}\tag{\theequation}}

\renewcommand{\S}{\mathbb{S}}

\makeatletter
\newtheorem*{rep@theorem}{\rep@title}
\newcommand{\newreptheorem}[2]{%
\newenvironment{rep#1}[1]{%
 \def\rep@title{#2 \ref{##1}}%
 \begin{rep@theorem}}%
 {\end{rep@theorem}}}
\makeatother
\newreptheorem{theorem}{Theorem}
\newreptheorem{lemma}{Lemma}
\newreptheorem{corollary}{Corollary}
\newreptheorem{proposition}{Proposition}

\usepackage{cleveref}
\newenvironment{delayedproof}[1]
 {\begin{proof}[\raisedtarget{#1}Proof of \Cref{#1}]}
 {\end{proof}}
\newcommand{\raisedtarget}[1]{%
  \raisebox{\fontcharht\font`P}[0pt][0pt]{\hypertarget{#1}{}}%
}
\newcommand{\proofref}[1]{\hyperlink{#1}{proof}}

\title{A Novel and Optimal Spectral Method for Permutation Synchronization}
\author{Duc Nguyen and Anderson Ye Zhang\\
~\\
University of Pennsylvania
}

\begin{document}
\maketitle
\begin{abstract} 
Permutation synchronization is an important problem in computer science that constitutes the key step of many computer vision tasks. The goal is to recover $n$ latent permutations from their noisy and incomplete pairwise measurements. In recent years, spectral methods have gained increasing popularity thanks to their simplicity and computational efficiency. 
Spectral methods utilize the leading eigenspace $U$ of the data matrix and its block submatrices $U_1,U_2,\ldots, U_n$ to recover the permutations. In this paper, we propose a novel and statistically optimal spectral algorithm. Unlike the existing methods which use $\{U_jU_1^\top\}_{j\geq 2}$, ours constructs an anchor matrix $M$ by aggregating useful information from all of the block submatrices and estimates the latent permutations through $\{U_jM^\top\}_{j\geq 1}$. This modification overcomes a crucial limitation of the existing methods caused by the repetitive use of $U_1$ and leads to an improved numerical performance.
To establish the optimality of the proposed method, we carry out a fine-grained spectral analysis and obtain a sharp exponential error bound that matches the minimax rate.
\end{abstract}

\section{Introduction}\label{sect:intro}
In permutation synchronization, the objective is to estimate $n$ latent permutations using noisy and potentially incomplete pairwise measurements among them. 
It is an important problem in computer vision and graphics where finding correspondence between sets of features across multiple images is a fundamental task with wide-ranging applications 
including image registration \cite{shen2002hammer}, shape matching \cite{berg2005shape}, multi-view matching \cite{gao2021isometric, maset2017practical}, detecting structures from motion \cite{agarwal2011building}, and so on. Various methods have been proposed for  permutation synchronization  including iterative algorithms \cite{gao2022iterative, zhou2015multi,chen2018projected}, semi-definite programming (SDP) \cite{huang2013consistent, chen2014near}, and  spectral methods \cite{pachauri2013solving,ling2022near,maset2017practical, shen2016normalized}. 
Compared to other approaches, spectral methods have gained increasing popularity and have been widely used in permutation synchronization thanks to their simplicity, fast computation speed, and impressive numerical performance. Despite the popularity, it remains unclear how well spectral methods perform theoretically and whether they achieve statistical optimality or not. In this paper, we address these questions by proposing a new and provably optimal spectral algorithm.

\paragraph{Problem Formulation.}
The permutation synchronization problem is formulated as follows.
Let $Z_1^*, \ldots, Z_n^* \in \Pi_d$ where $\Pi_d$ is the permutation group in $d$ dimension defined as:
\begin{equation}\label{eqn:Pi-def}
\Pi_d := \left\{ P \in \{0,1\}^{d\times d} \,:\, P^\top P = PP^\top = I_d  \right\}\,.
\end{equation}
We introduce missing and noisy data by assuming that, for each pair $1\leq j < k \leq n$, the observation $X_{jk} \in \R^{d\times d}$ satisfies
\begin{equation}\label{eqn:Xjk-def}
X_{jk} \define \begin{cases} Z_j^*(Z_k^*)^\top + \sigma W_{jk}, &\text{if } A_{jk} = 1,\\
0_{d\times d}, &\text{otherwise},
\end{cases}
\end{equation}
where $A_{jk} \in \{0,1\}$ independently follows $ \text{Bernoulli}(p)$ for some $p \in \left(0, 1\right]$; $\sigma \in \R^+$ controls the amount of noise; $W_{jk} \in \R^{d\times d}$ is a random matrix with each entry independently distributed according the standard normal distribution;  and $0_{d\times d}$ is the $d\times d$ matrix of all zeros. Roughly speaking, each $X_{jk}$ block, if not missing, is a noisy measurement of $Z_j^*(Z_k^*)^\top$ which is the `difference' between two permutation matrices. Denote $Z^*\define (Z_1^{*\top},\ldots, Z_n^{*\top})^\top\in \Pi_d^n$.  The goal is to estimate $Z^*$ given $\{ A_{jk} \}_{1\leq j < k \leq n}$ and $ \{X_{jk}\}_{1\leq j < k\leq n}$.
Note that $Z_1^*,\ldots, Z^*_n$ are identifiable only up to a global permutation. For any estimator $\hat Z =(\hat Z_1^\top, \ldots, \hat Z_j^\top)^\top\in \Pi_d^n$, its performance can be measured by the following normalized Hamming loss (modulo a global permutation transformation):
\begin{equation}\label{eqn:hamming}
\ell(\hat Z, Z^*) \define \min_{P \in \Pi_d} \,\frac{1}{n} \sum_{j=1}^n \indic{\hat Z_j \neq Z_j^* P^\top} \,.
\end{equation}
Note that model (\ref{eqn:Xjk-def}) has a matrix representation. %
The observation matrix $X \in \R^{nd\times nd}$ can be written  as 
\begin{equation}\label{eqn:X-def}
X = (A \otimes J_d) \circ Z^*(Z^*)^\top + \sigma (A\otimes J_d) \circ W,
\end{equation}
where $W \in \R^{nd\times nd}$ is a block-symmetric matrix with $W_{jj}:=0, W_{kj} := W_{jk}^\top,  \forall 1\leq j < k \leq n$; $A \in\{0,1\}^{n\times n}$ is a symmetric matrix with $A_{jj} := 0,A_{kj}:=A_{jk},\forall 1\leq j < k \leq n$; $J_d$ is the $d\times d$ matrix of all ones; $\otimes$ is the Kronecker product and $\circ$ is the Hadamard product.

\paragraph{Spectral Methods.}
Existing spectral methods \cite{pachauri2013solving,ling2022near,maset2017practical, shen2016normalized} use the eigendecomposition of $X$ followed by a rounding step to estimate the latent permutations. Let  $U=(U_1^\top,\ldots, U_n^\top)^\top \in \R^{nd\times d}$ be the matrix composing the top $d$ eigenvectors of $X$. That is, the columns of $U$ are the eigenvectors corresponding to the $d$ largest eigenvalues of $X$. Its $d\times d$ blocks are denoted as $U_1,\ldots, U_n$. The first block $U_1 \in \R^{d\times d}$ is then used as an `anchor' to obtain an estimator $\tilde Z_1, \ldots, \tilde Z_n \in \Pi_d$:
\begin{equation}\label{eqn:old-spectral}
\tilde Z_1 \define I_d, \quad \tilde Z_j \define \argmax_{P\in \Pi_d} \inner{P}{U_jU_1^\top} ,\forall j = 2,\ldots , n,
\end{equation}
where the optimization subproblem serves to round  $U_jU_1^\top$ into a permutation matrix and can be efficiently solved using the Kuhn-Munkres algorithm \cite{kuhn1955hungarian} (see Section \ref{sec:heuristic} for an intuitive description of the innerworking of this algorithm). To distinguish the existing algorithm from the one proposed in this paper, we refer to $\tilde Z:=(\tilde Z_1^\top,\ldots ,\tilde Z_n^\top)^\top$ as the `vanilla spectral estimator'. 
It is computationally efficient and has decent numerical performance. 

Despite of all the aforementioned advantages, the vanilla spectral method suffers 
from the repeated use of $U_1$ in constructing the estimator $\tilde Z_j$ for all $j\geq 2$. From a perturbation theoretical point of view, $U$ (along with its blocks $U_1,\ldots, U_n$) are approximations of their population counterparts. By using $U_jU_1^\top$, the estimation accuracy of $\tilde Z_j$ is determined by the approximation errors of \emph{both} $U_j$ and $U_1$. As a result, the approximation error of $U_1$ is carried forward in $\{\tilde Z_j\}_{j\geq 2}$ and deteriorates the overall numerical performance.

To overcome this crucial limitation of $\tilde Z$, we propose a new spectral method that avoids the use of $U_1$ as the anchor. Instead, we  construct an anchor matrix $M\in\mathr^{d\times d}$ by carefully aggregating useful information from all of $U$ and estimate the latent permutations by $\hat Z \define (\hat Z_1^\top, \ldots, \hat Z_n^\top)^\top$ where
\begin{align*}
\hat Z_j \define \argmax_{P\in \Pi_d} \inner{P}{U_jM^\top} ,\forall j\in[n].
\end{align*}
The construction of $M$ is built on an intuition that `averaging' information across all $n$ blocks of $U$  leads to a more accurate anchor with smaller variance than $U_1$. As a result, the estimation accuracy of $\hat Z_j$ is mostly determined by the approximation error of  $U_j$ only, which leads to an improved numerical performance. See Algorithm \ref{alg:main-algo} for the detailed implementation of the proposed method and Figure \ref{fig:err} for comparisons of numerical performances between the vanilla spectral estimator and ours.

\paragraph{Statistical Optimality.} By carrying out fine-grained spectral analysis, we establish a sharp upper bound for the theoretical performance  of the proposed method, summarized below in Theorem \ref{thm:intro}. We note that in this paper, $p,\sigma^2,d$ are not constants but  functions of $n$. This dependence can be more explicitly represented as $p_n,\sigma^2_n,d_n$. However, for simplicity of notation and readability, we choose to denote them as $p,\sigma^2,d$ throughout the paper. See Theorem \ref{thm:main_new} for a non-asymptotic and refined version where $d$ is also allowed to grow with $n$.
\begin{theorem}\label{thm:intro}
Assume $\frac{np}{\sigma^2}\rightarrow\infty$, $\frac{np}{\log^{3} n} \rightarrow\infty$ and $d = O(1)$. Then the proposed spectral method $\hat Z$ satisfies
\begin{align*}
\E \ell(\hat Z, Z^*) \leq \ebr{-(1-o(1)) \frac{np}{2\sigma^2}} + n^{-8}\,.
\end{align*}
\end{theorem}

The upper bound in Theorem \ref{thm:intro} consists of an exponential error term and a polynomial error term $n^{-8}$. Note that by properties of the normalized Hamming loss $\ell$, the polynomial error term is negligible. Considering this, our spectral algorithm achieves the minimax lower bound \cite{gao2022iterative} which states that if $\frac{np}{\sigma^2} \rightarrow \infty$, then $\inf_{Z'} \sup_{Z^*\in\Pi_d^n} \,\E \ell(Z', Z^*) \geq \ebr{-(1+o(1))\frac{np}{2\sigma^2}}$. This establishes the statistical optimality of the proposed method for the partial recovery of the latent permutations. Theorem \ref{thm:intro} immediately implies the threshold for exact recovery. When $np/(2 \sigma^2) >(1+\gamma)\log n$ for some constant $\gamma >0$, we have $\ell(\hat Z,Z^*)=0$ holds with high probability. According to the minimax lower bound, no estimator is able to recover $Z^*$ exactly with vanishing error if $np/(2 \sigma^2) <(1-\gamma)\log n$. As a result, simple but powerful, the proposed spectral method $\hat Z$ is an optimal procedure.

Theorem \ref{thm:intro}  allows the observations $\{X_{jk}\}$ to be missing at random as long as the probability $p$ satisfies $np\gg \log^3 n$. Note that in order to have a connected comparison graph $A$,  $np$ needs to be at least of order $\log n$. Compared to this condition, our assumption $np\gg \log^3 n$ requires some additional logarithm factor. The assumption $np \gg \sigma^2$ is the necessary and sufficient condition to achieve estimation consistency according to the minimax lower bound. Theorem \ref{thm:intro}  assumes that $d$, the dimension of each permutation matrix, is a constant.
To establish  Theorem \ref{thm:intro}, we first provide a block-wise $\ell_\infty$ perturbation analysis for all block submatrices $U_1,\ldots, U_n$, quantifying  the maximum deviation between them and their population counterparts (see Theorem \ref{thm:l_infty_new}). In addition, we give a theoretical justification for the usage of the anchor $M$ by showing that it achieves a negligible error (see Proposition \ref{prop:hat_H_new}). With both results, we investigate the tail behavior of each $U_jM^\top$ and eventually obtain the upper bound in  Theorem \ref{thm:intro}. We leverage the leave-one-out technique \cite{bean2013optimal, abbe2020entrywise} in our proofs.

\paragraph{Related Literature.} 

Permutation synchronization belongs to a broader class of group synchronization problems where the goal is to identify $n$ group objects  based on pairwise measurements among them. In recent years, spectral methods have been widely used and studied   in group synchronization problems. In \cite{zhang2022exact}, spectral methods are proved to be optimal for phase synchronization and orthogonal group synchronization in terms of squared $\ell_2$ losses. To obtain this result,  \cite{zhang2022exact} develops perturbation analysis toolkits to show that the leading eigenstructures can be well-approximated by its first-order approximation with a small $\ell_2$ error. However, the difference that permutations are discrete-valued while phases and orthogonal matrices are continuous is critical. For permutation synchronization, instead of $\ell_2$ perturbation analysis, we need to develop block-wise analysis in order to obtain sharp exponential rates. \cite{abbe2020entrywise} considers a $\mathbb{Z}_2$ synchronization problem where each object is $\pm 1$ and assumes that there is no missing data. It proves that a simple spectral procedure using signs of coordinates of the first eigenvector of the data matrix achieves the optimal threshold for the exact recovery of objects by $\ell_\infty$ analysis of the leading eigenvector. \cite{ling2022near} extends \cite{abbe2020entrywise}'s $\ell_\infty$ analysis to a permutation synchronization setting where there is no missing data and  each observation is corrupted with probability $q$ by a random permutation matrix. It shows that the vanilla spectral method achieves exact recovery  when $q$ satisfies certain conditions. Our work's novelty relative to \cite{ling2022near} is two-fold. Firstly, we propose a novel spectral algorithm that addresses a subtle but important limitation of the vanilla spectral algorithm, leading to a significant improvement in empirical performance. Secondly, our $\ell_\infty$ analysis is different from those  in \cite{abbe2020entrywise, ling2022near} as we need to consider the low signal-to-noise ratio  regime where  exact recovery is impossible but partial recovery is possible. We go beyond $\ell_\infty$ analysis and  study the tail behavior of each block of $U$ in order to obtain exponential error bounds for partial recovery. In addition, the presence of missing data in our model complicates the theoretical analysis as the magnitude and tail behavior of each block $U_j$ is not only related to the additive Gaussian noises but also the randomness of the Bernoulli random variables $\{A_{jk}\}_{k\neq j}$.

\paragraph{Notation.} For any positive integer $n$, define $[n] := \{1,2,\ldots, n\}$. Let $I_d$ denote the $d\times d$ identity matrix and $J_d$ denote the $d\times d$ matrix of all ones. Define $\mathcal{O}_d:=\{O\in\mathr^{d\times d}:OO^\top=O^\top O=I_d\}$ to be the set of all $d\times d$ orthogonal matrices.
Given $a, b\in \R$, let $a\lor b \define \max\{a, b\}$ and $a\land b \define \min\{a, b\}$. For a matrix $B \in \R^{d_1\times d_2}$, the Frobenius norm and the operator norm of $B$ are defined as $\fnorm{B} \define \left(\sum_{i=1}^{d_1}\sum_{j=1}^{d_2} B_{ij}^2 \right)^{1/2} $ and $\opnorm{B} \define \max_{u \in \S_{d_1}, v\in \S_{d_2}} u^\top B v $ where $\S_d \define \{v \in \R^d \,:\, \norm{v} = 1\}$ is the unit sphere in $d$ dimension and $\norm{.}$ is the Euclidean norm. For two matrices $A, B \in \R^{d_1\times d_2}$, let $\inner{A}{B} := \sum_{i\in [d_1]}\sum_{j\in [d_2]} A_{ij}B_{ij}$ denote the matrix inner product (this reduces to the usual vector inner product when $d_2 = 1$). For some $d_1, d_2 \in \mathbb{N}$, we use $\calN(\mu, \Sigma)$ to denote the normal distribution with mean $\mu \in \R^{d_1}$ and covariance $\Sigma \in \R^{d_1 \times d_1}$ and $\MN(S, \Sigma_1, \Sigma_2)$ to denote the matrix normal distribution \cite{de2004matrix} with mean parameter $S\in \R^{d_1\times d_1}$ and covariance parameters $\Sigma_1 \in \R^{d_1\times d_1}$, $\Sigma_2 \in \R^{d_2\times d_2}$. Define $(a)_+ \define  \max\{a, 0\}$. For the rest of the paper, we will use $e_i \in \{0,1\}^d$ to denote the vector with $1$ in the $i$-th entry and zero everywhere else. For any matrix $Y$ and any positive integer $i$, we denote $\lambda_i(Y)$ to be the $i$th largest eigenvalue of $Y$. 
For two positive sequences $\{a_n\}$ and $\{b_n\}$, $b_n\gtrsim a_n$ and $a_n=O(b_n)$ both mean $a_n\leq Cb_n$ for some constant $C>0$ independent of $n$. We also write $a_n=o(b_n)$ or $\frac{b_n}{a_n}\rightarrow\infty$ when $\limsup_n\frac{a_n}{b_n}=0$.  
Lastly, we use $\indic{.}$ to denote the indicator function.
 
\section{A Novel Spectral Algorithm}\label{sect:algo}
In this section, we first give a detailed implementation of the proposed method. In Section \ref{sec:heuristic}, we provide intuitions to explain why it works and how it improves upon the existing algorithm. In Section \ref{sec:numerical}, we compare the numerical performances of the proposed method with the vanilla spectral method on synthetic datasets. 

Algorithm \ref{alg:main-algo} describes our proposed spectral algorithm. The first step computes the top $d$ eigenvectors of the observation matrix $X$, which can be done using any off-the-shelf numerical eigendecomposition routine. The output of Step 1 is the empirical eigenspace $U \in \R^{nd\times d}$ corresponding to the top $d$ eigenvectors. Step 2 constructs the anchor $M$ used to recover the permutations by performing $d$-means clustering on the $nd$ rows of $U$  and extracting the estimated cluster centers. Step 3 recovers the underlying permutations via the Kuhn-Munkres algorithm. Note that (\ref{eqn:alg2}) in Step 3 can be interpreted as a projection of $U_j \hatH^\top$ onto $\Pi_d$. This is because all permutation matrices in $\Pi_d$ have the same Frobenius norm and (\ref{eqn:alg2}) can be equivalently written as
\begin{align}\label{eqn:alg3}
\hat Z_j = \argmin_{P\in\Pi_d}\fnorm{P - U_j M^\top}, \forall j\in[n].
\end{align}

\begin{algorithm}[ht]
\SetAlgoLined
\KwIn{Data matrix $X \in \mathr^{nd\times nd}$}
\KwOut{$n$ permutation matrices $\hat Z_1, \hat Z_2, \ldots, \hat Z_n\in\Pi_d$}
 \nl Obtain the top $d$ eigenvectors of $X$ as $U \in \mathr^{nd \times d}$\;
 \nl Denote the rows of $U$ as $v_1,\ldots,v_{nd}\in\mathr^{d\times 1}$ such that $U=(v_1^\top,\ldots, v_{nd}^\top)^\top$. Run $d$-means on $v_1,\ldots,v_{nd}$ and denote $\hat \mu_1,\ldots,\hat \mu_d\in\mathr^{d\times 1}$ as the cluster centers:
 \begin{align}
(\hat \mu_1,\ldots,\hat \mu_d)\define\argmin_{\mu_1,\ldots,\mu_d\in\mathr^{d\times 1}} \min_{z\in[d]^{nd}} \sum_{j=1}^{nd}\norm{v_j - \mu_{z_j}}^2.\label{eqn:alg1}
 \end{align}
Define $\hatH \define (\sqrt{n}\hat \mu_1^\top,\ldots,\sqrt{n}\hat \mu_d^\top)^\top\in\mathr^{d\times d}$ such that the rows of $\hatH$ are the $d$ centers multiplied by a $\sqrt{n}$ scaling\;
 \nl Compute for each $j\in[n]$,
 \begin{align}
 \hat Z_j \define \argmax_{P\in\Pi_d} \inner{P}{U_j \hatH^\top}. \label{eqn:alg2}
 \end{align}
\caption{A new spectral method for permutation synchronization. \label{alg:main-algo}}
\end{algorithm}

\subsection{Intuition}\label{sec:heuristic}
To understand Algorithm \ref{alg:main-algo}, we need to study $U$ through the lens of spectral perturbation theory. Denote
\begin{align}\label{eqn:U_star_def}
U^*\define (U^{*\top}_1,\ldots U^{*\top}_n)^\top \in\mathr^{nd\times n},\text{ where }U_j^* \define Z^*_j/\sqrt{n},\forall j\in[n].
\end{align}
Note that the expected value of the data matrix $X$ is equal to $\E X = p ((J_n - I_n) \otimes J_d) \circ Z^*(Z^*)^\top$. One could verify that $(\E X) U^* = (n-1)p U^*$. As a result, $U^*$ is the leading eigenspace of $\E X$ that includes its top $d$ eigenvectors. Since $X$ is a perturbed version of $\E X$, $U$ can also be seen as a perturbed version of $U^*$. However, this correspondence only holds with respect to an orthogonal transformation, as the leading eigenvalue of $\E X$, $(n-1)p $, has a multiplicity of $d$. That is, $UO^*\approx U^*$ for some $O^*\in\mathcal{O}_d$ (see Section \ref{sect:M} for the definition of $O^*$). As a result, the blocks $\{U_j\}$ are close to $\{U^*_j\}$ only up to a global orthogonal transformation, i.e., $U_jO^*\approx U^*_j$ for all $j\in[n]$.

To estimate the latent permutation $Z_j^*$, the vanilla spectral method uses the product $U_jU_1^\top$ as $U_jU_1^\top = (U_jO^*)(U_1O^*)^\top \approx U^*_j(U_1^*)^\top = Z^*_j(Z^*_1)^\top/n$. Intuitively, if $U_jU_1^\top$ is very close to  $Z^*_jZ^*_1/n$ then by (\ref{eqn:old-spectral}), we have $\tilde Z_j = Z^*_j(Z^*_1)^\top$. As a result, when the perturbation between $UO^*$ and $U^*$ is sufficiently small, the vanilla spectral method is able to recover $Z^*$ up to a global permutation $Z^*_1$. However, as hinted at before, the use of the product $\{U_jU_1^\top\}$ in the vanilla spectral method inevitably leads to a crucial and fundamental limitation. For each $j\in[n]$, write $U_jO^*= U_j^*+\xi_j$ such that $\xi_j$ can be interpreted as the approximation `noise' of $U_jO^*$ with respect to $U_j^*$. Then
\begin{align*}
   U_jU_1^\top =   (U_jO^*)(U_1O^*)^\top = (U^*_j + \xi_j)(U^*_1 + \xi_1)^\top= U_j^*(U_1^*)^\top + \xi_j(U^*_1)^\top + U^*_j  \xi_1^\top + \xi_j \xi_1^\top.
\end{align*}
That is, $ U_jU_1^\top$ is an approximation of $U_j^*(U_1^*)^\top$ with an error $\xi_j(U^*_1)^\top + U^*_j  \xi_1^\top$ (the higher-order  term $\xi_j \xi_1^\top$ is ignored).  As a function of $   U_jU_1^\top$, the estimation accuracy of $\tilde Z_j$ is determined by both $\xi_j$ and $\xi_1$. As a result, the error caused by $\xi_1$  is carried forward in $\{\tilde Z_j\}_{j\geq 2}$ and impairs the numerical performance of the overall algorithm (see Figure \ref{fig:err}). 
Additionally, using $U_1$ as the anchor makes the performance of the vanilla spectral algorithm less stable because it highly depends on the accuracy of $U_1$.

Algorithm \ref{alg:main-algo} overcomes this crucial limitation of the vanilla spectral method by constructing an anchor that `averages' information across all rows of $U$ instead of just using its first block $U_1$. The key insight into our construction is to recognize the special structure of the permutation synchronization problem where the true eigenspace $U^*$ consists of $d$ unique rows $e_1^\top/\sqrt n,\ldots, e_d^\top/\sqrt n$, each of cardinality exactly $n$. 
The empirical eigenspace $U$, being a noisy estimate of $U^* {O^*}^\top$, thus exhibits clustering structures where the cluster centers are the transformed rows $\{e^\top_i {O^*}^\top/\sqrt n \}_{i=1}^d$
Clustering algorithms such as the $d$-means algorithm can be used to estimate the cluster centers, and thereby the rotation matrix $O^*$, accurately. With $M$ being an accurate approximation of ${O^*}^\top$ (up to a row permutation), we have
\begin{align*}
 U_jM^\top \approx U_jO^* = U^*_j + \xi_j,
\end{align*}
and consequently the estimation accuracy of $\hat Z_j$ is only related to $\xi_j$ where as the accuracy of $U_jU_1^\top$ depends on both $\xi_j$ and $\xi_1$. The reduced noise in $U_jM^\top$ leads to an improved and more stable numerical performance.

\subsection{Numerical Analysis}\label{sec:numerical}
To verify our intuitions and showcase the improved performance of our spectral method over the vanilla spectral method, we perform experiments on synthetic data following the model in (\ref{eqn:Xjk-def}). 

\begin{figure}[!ht]
\begin{subfigure}[t]{0.45\textwidth}
\centering
   \includegraphics[width=1\linewidth]{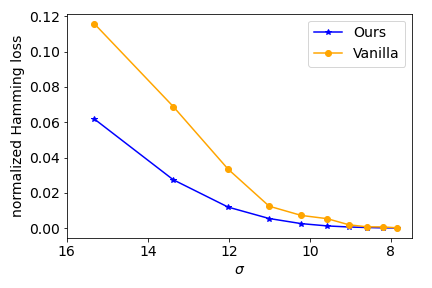}
   \caption{
   \label{fig:hamming}}
\end{subfigure}
\hspace{0.5cm}
\begin{subfigure}[t]{0.45\textwidth}
\centering
   \includegraphics[width=1\linewidth]{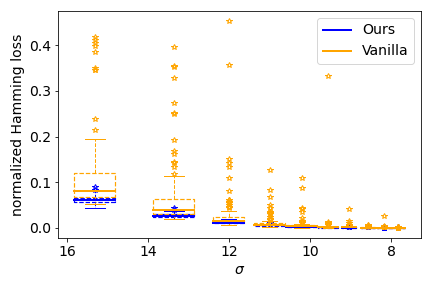}
   \caption{
   \label{fig:box-plot}}
\end{subfigure}\\

\begin{subfigure}[t]{0.45\textwidth}
\centering
   \includegraphics[width=1\linewidth]{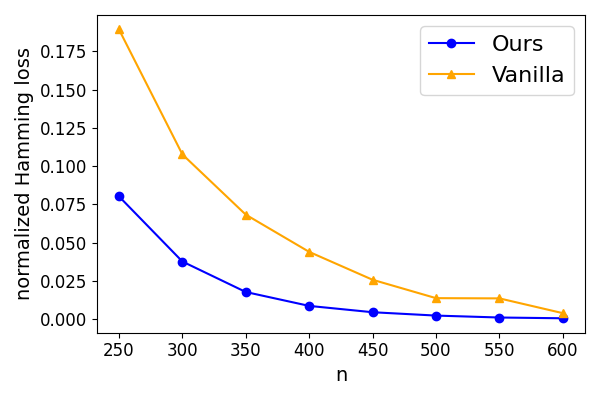}
   \caption{
   \label{fig:loss-vs-n}}
\end{subfigure}
\hspace{0.5cm}
\begin{subfigure}[t]{0.45\textwidth}
\centering
   \includegraphics[width=1\linewidth]{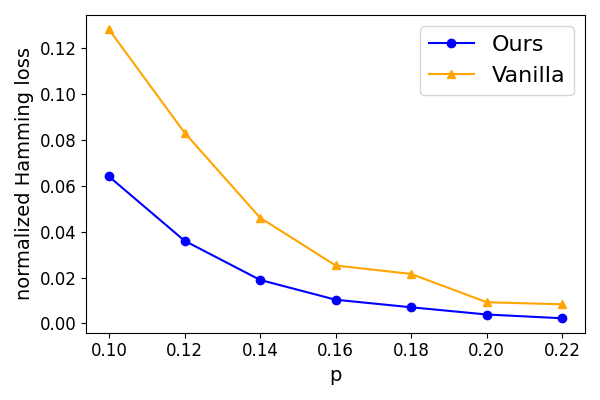}
   \caption{
   \label{fig:loss-vs-p}}
\end{subfigure}

\caption{{Comparisons between our method and the vanilla spectra method on synthetic data.} \label{fig:err}}
\end{figure}

In Figures \ref{fig:hamming}-\ref{fig:loss-vs-p}, we compare the performance of Algorithm \ref{alg:main-algo} against the vanilla spectral algorithm. In each figure, we vary a single parameter while keeping  all of the other parameters fixed: in Figures \ref{fig:hamming}-\ref{fig:box-plot}, we set $d = 2, p = 0.5, n = 2048$ and vary $\sigma$; in Figure \ref{fig:loss-vs-n}, we set $\sigma = 1, p = 0.1, d = 3$ and vary $n$; in Figure \ref{fig:loss-vs-p}, we set  $n = 200, \sigma = 2, d = 2$ and vary $p$. In Figures \ref{fig:hamming} and \ref{fig:loss-vs-n}-\ref{fig:loss-vs-p}, each line shows the average over 100 independent trials. One can see that our spectral method outperforms the vanilla spectral method across all settings.

The box plot in Figure \ref{fig:box-plot} shows the distribution of the losses across 100 trials for both methods. The boxes extend from the first quartile to the third quartile of the observed losses with colored bolded lines at the medians and dotted points as the outliers. It is clear that the performance of the vanilla algorithm is highly variable across trials. On the contrary, our method has a smaller variance and is more precise. Figure \ref{fig:err} reflects our intuition that because the vanilla algorithm repeatedly uses $U_1$ to estimate the latent permutations, its performance is highly dependent on the approximation error of a single block $U_1$. On the other hand, because our algorithm constructs an anchor that averages information across all $n$ blocks, it is more stable and consistently outperforms the vanilla algorithm with a much tighter error spread.

\section{Theoretical Guarantees}\label{sect:analysis}
In this section, we  establish theoretical guarantees for Algorithm \ref{alg:main-algo}. Similarly to previous analysis \cite{ling2022near},  without loss of generality\footnote{This is because for any $X$  as defined in (\ref{eqn:Xjk-def}) with an arbitrary $Z^*\in\Pi_d^n$, it can be transformed into $X'\in\mathr^{nd\times nd}$ by letting $X'_{jk} := (Z^*_j)^\top X_{jk} Z^*_k,\forall1 \leq j,k\leq n$. 
One can show that implementing Algorithm 1 on either $X$ or $X'$ yields identical results, up to a global permutation. In addition, it can be shown that $X'$ also adheres to the definition in (\ref{eqn:Xjk-def}), where all latent permutation matrices are the identity matrix $I_d$.
}, we assume that $Z^*_j = I_d$ for all $j\in[n]$.
That is, all the latent permutations are the identity matrix. In this way, the population eigenspace $U^*$ has a simpler expression with
\begin{align}\label{eqn:U_star_simplified}
U^*_j =I_d/\sqrt{n},\forall j\in[n].
\end{align}

We first give a justification for the choice of the anchor $M$, the key algorithmic novelty of the proposed algorithm. We then establish the statistical optimality of the proposed method in Section \ref{sec:optimality}.  As a preview of our proofs, 
we then present an intuitive analysis of the block-wise perturbation analysis for $U$ which will serve as the starting point of our detailed proofs. Last, we propose a practical modification of Algorithm \ref{alg:main-algo} in Section \ref{sect:theory-approx} which uses approximate clustering instead of exact clustering, allowing the algorithm to remain scalable when $d$ diverges. The modified spectral algorithm enjoys similar theoretical guarantee as Algorithm \ref{alg:main-algo} under mild additional assumptions.

\subsection{Accuracy of the Anchor $M$}\label{sect:M}
At a high level, the effectiveness of $M$ as the anchor stems from the property that $UM^\top$ is close to $U^*$ up to a global permutation. However, this is not obvious at first glance.

Note that there is a close relationship between the empirical eigenspace $U$ and the population eigenspace $U^*$, modulo an appropriate orthogonal transformation. Readers familiar with spectral analysis literature might anticipate the natural choice for this transformation to be the sign matrix \cite{chen2021spectral}, defined as $O^* := \argmin_{O \in \mathcal{O}_d} \fnorm{UO - U^*}$. However, our analysis employs a different approach, using the following matrix \cite{abbe2020entrywise}:
\begin{equation}
H \define U^\top U^*\,.
\end{equation}
This choice is strategic; although $H$ is not orthonormal, it offers a straightforward and explicit expression that facilitates easier manipulation than $O^*$. Additionally, $O^*$, by definition, represents the orthogonal matrix that most closely approximates $H$, and the distance between them is small. We first show (see (\ref{eqn:5}) of Lemma \ref{lem:event-E-prime_new}):
\begin{align}\label{eqn:new_1}
\opnorm{UH-U^*}  =O\br{\frac{{1+\sigma\sqrt{d}}}{\sqrt{np}}}
\end{align}
holds with high probability.
Then following the intuition given in Section \ref{sec:heuristic} and applying state-of-the-art clustering analysis, we have the following proposition which states that the anchor $M$ is in fact an approximation of $H^\top$ up to some permutation matrix $\hat P$ where $\hat P \define \argmin_{P\in \Pi_d} \fnorm{M - PH^\top}$.
\begin{proposition}\label{prop:hat_H_new}
There exist constants $C_1,C_2,C_3>0$ such that if $\frac{np}{\log n}>C_1$ and $\frac{np}{(\sqrt{d} + \sigma d)^2}>C_2$ then the anchor $M$ constructed in Algorithm \ref{alg:main-algo} satisfies
  \begin{align}\label{eqn:21}
    \min_{P\in\Pi_d}\fnorm{\hatH - PH^\top}\leq  \frac{C_3(\sqrt{d} + \sigma d)}{\sqrt{np}},
\end{align}
with probability at least $1- n^{-10}$.
\end{proposition}
The existence of $\hat P$ is due to the fact that the clusters are identifiable only up to a permutation in cluster analysis. More importantly, the exact value $\hat{P}$ does not affect the final error bound in (\ref{eqn:hamming}) because the Hamming loss is defined with respect to an unknown optimal global permutation. (\ref{eqn:21}) shows that the estimation error goes to 0 when $\frac{np}{(\sqrt{d}+\sigma d)^2}$ grows. As a comparison, the operator norm of the matrix $H$ is `of a constant order' as it is close to an orthogonal matrix (see (\ref{eqn:6}) of Lemma \ref{lem:event-E-prime_new}). Hence, the error incurred by $M$ as an estimate of $H^\top$ is of diminishing proportion. It is also worth mentioning that Proposition \ref{prop:hat_H_new} holds under weaker assumptions compared to Theorem \ref{thm:intro}. It only requires $np \gtrsim \log n$ and allows $d$ to grow as long as $np \gtrsim (\sqrt{d} + \sigma d)^2$. With Proposition \ref{prop:hat_H_new}, we can have a decomposition of the quantity $U_jM^\top$, the key quantity in (\ref{eqn:alg2}).
\begin{align*}
U_jM^\top -U_j^* \hat P^\top &= U_j(\hat PH^\top + M - \hat PH^\top)^\top -U_j^* \hat P^\top \\
&= \br{U_j H - U_j^*}\hat P^\top + U_j(M - \hat PH^\top)^\top\\
&\approx  \br{U_j H - U_j^*}\hat P^\top. \numberthis \label{eqn:approx1}
\end{align*}
The approximation (\ref{eqn:approx1}) is due to (\ref{eqn:21}) as the term $ U_j(M - \hat PH^\top)^\top$ turns out to be negligible. Hence, the difference $U_jM^\top -U_j^* \hat P^\top $ is primarily about the perturbation $U_j H - U_j^*$. To analyze our spectral method, we need to have a deep understanding on the behaviors of  the blockwise perturbation $U_j H - U_j^*$.

\subsection{Statistical Optimality}\label{sec:optimality}

In Theorem \ref{thm:l_infty_new}, we first derive a block-wise $\ell_\infty$ upper bound for $UH-U^*$, i.e.,  an upper bound for $\opnormt{U_j H - U^*_j}$ that holds uniformly across all $j\in[n]$.

\begin{theorem}\label{thm:l_infty_new}
There exist constants $C,C',C''>0$ such that if $\frac{np}{(\sigma^{4/3}  \vee 1)\log n} >C$ and $\frac{np}{(\sqrt{d} + \sigma d)^2}\geq C'$, we have
\begin{align}\label{eqn:22}
\max_{j\in[n]} \opnorm{U_j H - U^*_j}\leq \frac{ C''}{\sqrt{n}}\br{\sqrt{\frac{\log n }{np}} + \frac{\sigma \sqrt{d}}{\sqrt{np}} + \frac{\sigma \sqrt{\log n}}{\sqrt{np}}}
\end{align}
with probability at least $1-6dn^{-9}$.
\end{theorem}

The upper bound in (\ref{eqn:22}) is equal to the upper bound of $\opnorm{UH - U^*}$ in (\ref{eqn:new_1}) multiplied by a $\sqrt{(\log n)/n}$ factor. This is because $UH-U^*$ consists of $n$ blocks $\{U_jH-U^*_j\}_{j\in[n]}$  and all blocks behave similarly.  The magnitude of each block is  then on average $1/\sqrt{n}$ of that of the whole matrix, and the $\sqrt{\log n}$ factor is due to the use of union bound to control the supremum. Theorem \ref{thm:l_infty_new} assumes $\frac{np}{(\sigma^{4/3}  \vee 1)\log n} \gtrsim 1$. A sufficient condition is ${np}/{\log^3 n}\gtrsim 1$  as $np/\sigma^2\gtrsim 1$ is implied by the other assumption  ${np}/{(\sqrt{d} + \sigma d)^2}\gtrsim 1$.

With the upper bound for each $\opnormt{U_j H - U^*_j}$ derived, we further study the tail behavior of each difference $U_j H - U^*_j$. This leads to a sharp theoretical analysis of the performance of the proposed spectral method. The main theoretical result of this paper is stated below in Theorem \ref{thm:main_new}.
\begin{theorem}\label{thm:main_new} 
There exist constants $\bar C',\bar C'',\bar C''', \bar C''''>0$ such that if $\frac{np}{d^2\log n}>\bar C'$, $\frac{np}{\sigma^2d^3}>\bar C''$ and $\frac{np}{\log^3n} > \bar C''''$ then the estimates $\hat Z_1,\ldots, \hat Z_n$ of Algorithm \ref{alg:main-algo} satisfy
\begin{align*}
\E \ell(\hat Z,Z^*)\leq \ebr{-\br{1-\bar C''''\br{\br{\frac{\sigma^2d^3}{np}}^\frac{1}{4} + \br{\frac{d^2\log n}{np}}^\frac{1}{4} }}\frac{np}{2\sigma^2}}  + n^{-8}.
\end{align*}
\end{theorem}

Theorem \ref{thm:main_new} is a non-asymptotic and refined version of Theorem \ref{thm:intro} stated in the introduction.
By letting $np/(\sigma^2d^3)$ and $np/(d^2\log n)$ go to infinity, the exponential term in Theorem \ref{thm:main_new} takes an asymptotic form of $\ebr{-(1-o(1))\frac{np}{2\sigma^2}}$ and matches with the minimax lower bound. In this way, we establish the statistical optimality of the proposed spectral method. 
In Theorem \ref{thm:main_new}, $d$ is allowed to grow with $n$ but not too fast. 
\cite{gao2022iterative} considers the full observation setting where $p = 1$ and studies an EM-like iterative algorithm, requiring that $n/(\sigma^2 d \log d) \rightarrow\infty$. Our conditions are more stringent, reflecting the added complexity introduced by missing data, which necessitates stronger assumptions to apply concentration inequalities effectively.

While the proof of Theorem \ref{thm:main_new} is complicated, one can still develop an intuition as to why we can achieve the error bound $\ebr{-(1-o(1))\frac{np}{2\sigma^2}}$. 
At a high level, we will prove and use a variant of the linear approximation \cite{abbe2020entrywise} $UH - U^* \approx XU^*\Lambda^{-1} - U^*$ where $\Lambda$ is the diagonal matrix of the leading eigenvalues of $X$ (see (\ref{eqn:decomposition1})). A further decomposition using the structure of $X$ reveals that 
\begin{align}\label{eqn:approx2}
U_jH-U_j^*\approx  \frac{\sigma}{\sqrt{n}} \,\sum_{k\neq j} A_{jk} W_{jk} \Lambda^{-1}  \approx  \frac{1}{\sqrt{n}} \frac{\sigma}{np} \sum_{k\neq j}A_{jk}W_{jk}.
\end{align}
Recall that $U_j^*= I_d/\sqrt{n}$ according to (\ref{eqn:U_star_simplified}). Then $\sqrt{n }U_jH - I_d \approx \frac{\sigma}{np} \sum_{k\neq j}A_{jk}W_{jk}$ which follows a Gaussian distribution conditioned on $\{A_{jk}\}_{k\neq j}$. Since $\sum_{k\neq j}A_{jk}$ concentrates around $(n-1)p$, roughly speaking, $\sqrt{n }U_jH - I_d$ is a random matrix with each entry i.i.d. following $\mathcal{N}(0,\frac{\sigma^2}{np})$. The Gaussian tail is used to characterize the probability of the event when $\sqrt{n }U_jH$ considerably deviates away from $I_d$ and eventually leads to a probability bound of $\ebr{-(1-o(1))\frac{np}{2\sigma^2}}$.

\subsection{Block-wise Decomposition of $UH-U^*$}
In this section, we provide a block-wise decomposition of $UH-U^*$. The decomposition is the key towards the block-wise analysis in Section \ref{sec:optimality} and provides insights on how Theorem \ref{thm:l_infty_new} and Theorem \ref{thm:main_new} are established. Recall that we let $Z^*_j = I_d$ for all $j\in[n]$. Then, (\ref{eqn:X-def}) becomes
\begin{align*}
X = (A \otimes I_d) + \sigma (A \otimes J_d) \circ W.
\end{align*}
For each $j\in[n]$, define $X_j\define (X_{j1},\ldots, X_{jn})\in\mathr^{d\times nd}$ to be the $j$-th block row of $X$ and define $W_j \in \mathr^{d\times nd}$ analogously for $W$. Then we have
\begin{align}\label{eqn:1}
X_j = (A_j \otimes I_d) + \sigma (A_j \otimes J_d) \circ W_j.
\end{align}
As a result, $X =(X_1^\top,\ldots, X_n^\top)^\top$.  Define $\Lambda\in\mathr^{d\times d}$ to be the diagonal matrix of the leading eigenvalues of $X$. That is,
$\Lambda_{ii}\define \lambda_i(X)$ and $\Lambda_{ik}\define 0$ for all $1\leq i\neq k\leq d$.
Then we have
\begin{align*}
UH\Lambda - XU^*  & = U(H\Lambda - \Lambda H) + U\Lambda H -XU^*\\
& = U(H\Lambda - \Lambda H) + XU H -XU^*\\
&= U(H\Lambda - \Lambda H) + X(U H - U^*).
\end{align*}
Multiplying both sides by $\Lambda^{-1}$ and rearranging the terms, we have
\begin{align}\label{eqn:decomposition1}
UH - U^* = U(H\Lambda - \Lambda H) \Lambda^{-1}  + X(U H - U^*) \Lambda^{-1} + XU^*\Lambda^{-1} - U^* .
\end{align}
The above display involves $ X(U H - U^*) $ where $X$ and $UH-U^*$ are dependent on each other. To decouple the dependence, we approximate $UH-U^*$ by its leave-one-out counterparts.
Consider any $j\in[n]$. Define $X^{(j)}\in\mathr^{nd\times nd}$ such that
\begin{align*}
X^{(j)}_{ik} \define\quad \begin{cases}
X_{ik}, &\forall i,k\neq j,\\
0_{d\times d},  &\text{otherwise.}
\end{cases}
\end{align*}
In addition, let $U^{(j)}\in\mathr^{nd\times d}$ be the matrix including the leading $d$ eigenvectors of $X^{(j)}$. As a consequence, $X^{(j)}, U^{(j)}$ are independent of $\{A_{jk}\}_{k\neq j}$ and $\{W_{jk}\}_{k\neq j}$. 
Define
$
H^{(j)}\define U^{(j)\top}U^*\in\mathr^{d\times d}.
$
Then
\begin{align*}
UH-U^* &= UH - U^{(j)}H^{(j)} + U^{(j)}H^{(j)} - U^*  = (UU^\top - U^{(j)}U^{(j)\top})U^*+ U^{(j)}H^{(j)} - U^* .
\end{align*}
After plugging it into the right-hand side of (\ref{eqn:decomposition1}), we have
\begin{align*}
UH - U^* & = U(H\Lambda - \Lambda H) \Lambda^{-1}  + X(UU^\top - U^{(j)}U^{(j)\top})U^* \Lambda^{-1} \\
&\quad  +  X(U^{(j)}H^{(j)} - U^*) \Lambda^{-1}  + XU^*\Lambda^{-1} - U^*.
\end{align*}
Then, the $j$th block matrix of $UH - U^*$ satisfies
\begin{align}
U_jH - U^*_j & = \underbrace{U_j(H\Lambda - \Lambda H) \Lambda^{-1}  + X_j(UU^\top - U^{(j)}U^{(j)\top})U^* \Lambda^{-1}}_{\definei B_j} \nonumber\\
&\quad  +  X_j(U^{(j)}H^{(j)} - U^*) \Lambda^{-1}  + X_jU^*\Lambda^{-1} - U^*_j.\label{eqn:decomposition2}
\end{align}
The last two terms in (\ref{eqn:decomposition2}) can be further decomposed. 
Using (\ref{eqn:1}), we have
\begin{align*}
X_j(U^{(j)}H^{(j)} - U^*) \Lambda^{-1}&= \br{(A_j \otimes I_d) + \sigma (A_j \otimes J_d) \circ W_j}(U^{(j)}H^{(j)} - U^*) \Lambda^{-1}\\
&=  \underbrace{\sum_{k\neq j} A_{jk} (U^{(j)}_kH^{(j)} - U^*_k)\Lambda^{-1}}_{\definei F_{j1}}  +   \underbrace{\sigma \sum_{k\neq j} A_{jk} W_{jk}(U^{(j)}_kH^{(j)} - U^*_k)\Lambda^{-1}}_{\definei F_{j2}}  ,
\end{align*}
and
\begin{align*}
X_jU^*\Lambda^{-1} - U^*_j &= \br{(A_j \otimes I_d) + \sigma (A_j \otimes J_d) \circ W_j}U^*\Lambda^{-1} - U^*_j \\
&= \left(\sum_{k\neq j} A_{jk} U^*_k\right)\Lambda^{-1}  - U_j^*+ \sigma \,\sum_{k\neq j} A_{jk} W_{jk} U^*_k\Lambda^{-1} \\
&= \underbrace{\frac{1}{\sqrt{n}}\br{\left(\sum_{k\neq j} A_{jk} \right)\Lambda^{-1}  - I_d}}_{\definei G_{j1}} + \underbrace{ \frac{\sigma}{\sqrt{n}} \,\sum_{k\neq j} A_{jk} W_{jk} \Lambda^{-1} }_{\definei G_{j2}} ,
\end{align*}
where the last equation is due to Lemma \ref{lem:population}.
As a result, we have a decomposition of $U_jH - U^*_j$
\begin{align}\label{eqn:decomposition3}
U_jH - U^*_j & = 
B_j  + F_{j1} + F_{j2} + G_{j1} + G_{j2},
\end{align}
holds for all $j\in[n]$.

The block-wise decomposition (\ref{eqn:decomposition3}) of $UH-U^*$ is the starting point to establishing both Theorem \ref{thm:l_infty_new} and Theorem \ref{thm:main_new}. 
Note that we have a mutual independence among $\{A_{jk}\}_{k\neq i}$, $\{W_{jk}\}_{k\neq i}$, and $U^{(j)}H^{(j)} -U^*$ in the definitions of $F_{j1}$ and $F_{j2}$, which is crucial to obtaining sharp bounds and tail probabilities for them. Theorem \ref{thm:l_infty_new} is proved by establishing upper bounds for the operator norm of $B_j, F_{j1}, F_{j2}, G_{j1} $, and $G_{j2}$. By further analyzing their tail bounds, we establish Theorem \ref{thm:main_new}. Among these terms, $G_{j2}$ is the one contributing to the exponential error bound in Theorem \ref{thm:main_new}, as we illustrate in (\ref{eqn:approx2}).

\subsection{Approximate Clustering}\label{sect:theory-approx}
Algorithm \ref{alg:main-algo} involves a $d$-means clustering (\ref{eqn:alg1}).
There exists an $O(n^{(O(d^2))})$ algorithm using a weighted Voronoi diagram \citep{inaba1994applications} that finds the globally optimal solution. In the case where $d$ is a constant, this algorithm enjoys polynomial time complexity. However, this approach quickly becomes impractical as $d$ grows. In fact, solving the general $d$-means optimization problem exactly is NP-hard \citep{aloise2009np,mahajan2012planar}. 
A practical solution is to use an approximate algorithm which guarantees a $(1+\epsilon)$-optimal solution to (\ref{eqn:alg1}) under a polynomial time constraint. That is, let $\check \mu_1, \ldots, \check\mu_d \in\mathr^{d\times 1}$ and $\check z \in [d]^{nd}$ denote a $(1+\epsilon)$-approximate solution. 
That is, they satisfy 
\begin{equation}\label{eqn:approx-cls-err}
  \sum_{j=1}^{nd}\norm{v_j - \check\mu_{\check z_j}}^2 \leq (1+\epsilon)\min_{\mu_1,\ldots,\mu_d\in\mathr^{d\times 1}} \min_{z\in[d]^{nd}} \sum_{j=1}^{nd}\norm{v_j - \mu_{z_j}}^2,
\end{equation}
Define an approximate anchor as 
\begin{equation}\label{eqn:approx-M}
 \check M \define \left(\sqrt n\check \mu_1^\top,\ldots, \sqrt n\check \mu_d^\top \right)^\top \in\mathr^{d\times d}.
\end{equation} 
We will refer to the anchor $M$ used in Algorithm \ref{alg:main-algo} as the exact anchor.
Using the randomized $d$-means algorithm of \citep{kumar2010linear}, we obtain a $(1+\epsilon)$-approximate solution to (\ref{eqn:alg1}) in $O(nd2^{(d/\epsilon)^{O(1)}})$ time, which can be considerably faster than the exact clustering approach. We also emphasize that the third step of Algorithm \ref{alg:main-algo} is agnostic of the anchor. Therefore, we simply replace $M$ by $\check M$ with no change to the remaining of the algorithm, as shown in the following.
\begin{algorithm}[ht]
\SetAlgoLined
\KwIn{Data matrix $X \in \mathr^{nd\times nd}$, approximation factor $\epsilon>0$ }
\KwOut{$n$ permutation matrices $\check Z_1, \check Z_2, \ldots, \check Z_n\in\Pi_d$}
 \nl Obtain the top $d$ eigenvectors of $X$ as $U \in \mathr^{nd \times d}$\;
 \nl Denote the rows of $U$ as $v_1,\ldots,v_{nd}\in\mathr^{d\times 1}$. Run $(1+\epsilon)$-approximate $d$-means of clustering on $v_1,\ldots,v_{nd}$ and denote $\check{\mu}_1,\ldots, \check{\mu}_d \in\mathr^{d\times 1}$ as the approximate cluster centers which satisfy (\ref{eqn:approx-cls-err})\;
\nl Compute for each $j\in[n]$,
\begin{align}
\check Z_j \define \argmax_{P\in\Pi_d} \inner{P}{U_j \check\hatH^\top},
\end{align}
where $\check M$ is defined per (\ref{eqn:approx-M}).
\caption{Spectral method with approximate clustering. \label{alg:approx-algo}}
\end{algorithm}

With a small additional assumption regarding $(1+\epsilon)$, in Theorem \ref{thm:main_approx}, we show that the theoretical guarantee of the modified spectral algorithm is the same as that of Algorithm \ref{alg:main-algo}. The additional assumption $ \frac{np}{\sigma^2 d (1+\epsilon)}\gtrsim 1$ is due to the use of the approximate anchor $\check M$ instead of the exact one $M$. The proof of  Theorem \ref{thm:main_approx} is nearly identical to that of Theorem \ref{thm:main_new} with minor modifications.

\begin{theorem}\label{thm:main_approx} 
There exist constants $\check C_1,\check C_2,\check C_3, \check C_4, \check C_5>0$ such that if $\frac{np}{d^2\log n}>\check C_1$, $\frac{np}{\sigma^2d^3}>\check C_2$, $\frac{np}{\log^3n} > \check C_3$, and $\frac{np}{\sigma^2 d} \geq \check C_4(1+\epsilon)$, then the estimates  $\check Z_1,\ldots, \check Z_n$ of Algorithm \ref{alg:approx-algo} satisfy
\begin{align*}
\E \ell(\check Z,Z^*)\leq \ebr{-\br{1-\check C_5\br{\br{\frac{\sigma^2d^3}{np}}^\frac{1}{4} + \br{\frac{d^2\log n}{np}}^\frac{1}{4} }}\frac{np}{2\sigma^2}}  + n^{-8}.
\end{align*}
\end{theorem}

\section{Proofs}

\subsection{Preliminaries and Useful Inequalities}\label{sect:preliminaries}

We first present a lemma which enumerates a number of useful intermediate inequalities. Conditioned on certain inequalities, the statements of many subsequent lemmas and theorems hold deterministically. This allows us to present the proofs of the lemmas and theorems in more intuitive ways. The proof of Lemma \ref{lem:event-E-prime_new} is deferred to a later section.

\begin{lemma}\label{lem:event-E-prime_new} 
There exist constants $C_0'>0, C_0 > 7$ such that if $\frac{np}{\log n} > C_0'$, then the following event holds
\begin{align*}
 \calE &= \left\{ \opnorm{A - \E A} \leq C_0\sqrt{np} ,\;   \opnorm{(A \otimes J_d)\circ W} \leq C_0\sqrt{npd},\; \max_{j\in[n]}\abs{\sum_{k\neq j}A_{jk}-np}\leq C_0\sqrt{np\log n}  \right\}\,
\end{align*}
with probability at least $1- n^{-10}$. Under the event $\calE$, we have
\begin{align}
\opnorm{X - (\E A \otimes I_d)} ,\max_{i\in[d]} \abs{\lambda_i(X) - (n-1)p} ,\max_{d+1\leq i \leq nd} \abs{\lambda_i(X)+p} &\leq  C_0\br{1+\sigma\sqrt{d}}\sqrt{np},\label{eqn:2}\\
\max_{j\in[n]} \opnorm{X_j} \leq p\sqrt{n} + C_0\br{1+\sigma\sqrt{d}}\sqrt{np}.\label{eqn:7}
\end{align}
If $\frac{np}{(1+\sigma \sqrt{d})^2 }\geq 64C_0^2$ is further assumed, under the event $\calE$, the following hold for $\Lambda, U$, and $H$:
\begin{align}
\opnorm{\Lambda} =\lambda_1(X) &\leq \frac{9np}{8},\label{eqn:8}\\
\opnorm{\Lambda^{-1}}^{-1} =\lambda_d(X) &\geq \frac{7np}{8},\label{eqn:3}\\
\lambda_d(X)- \lambda_{d+1}(X) &\geq \frac{3np}{4},\label{eqn:4}\\
 \min_{O\in \mathcal{O}_d}\opnorm{H - O}, \opnorm{UH-U^*} ,  \opnorm{UU^\top - U^* U^{*\top} } &\leq   \frac{8C_0\br{1+\sigma\sqrt{d}}}{7\sqrt{np}},\label{eqn:5}\\
 \opnorm{H\Lambda - \Lambda H} &\leq  2C_0\br{1+\sigma\sqrt{d}}\sqrt{np},\\
\opnorm{H^{-1}} &\leq \frac{4}{3},\label{eqn:6}
\end{align}
 the following hold any $j\in[n]$:
\begin{align}
\opnorm{G_{j1}} &\leq \frac{2C_0}{\sqrt{n}}\br{\sqrt{\frac{\log n}{np}} + \frac{\sigma \sqrt{d}}{\sqrt{np}}},\label{eqn:16}\\
 \opnorm{U_j} &\leq \frac{4}{3}\br{\opnorm{U_jH-U^*_j}+\frac{1}{\sqrt{n}}},\label{eqn:10} \\
 \opnorm{\loo{U}(\loo{U})^\top -UU^\top} &\leq 6\br{\opnorm{U_jH -U^*_j}+\frac{1}{\sqrt{n}}}, \label{eqn:14}\\
\opnorm{\loo{U}\loo{H} -U^*}&\leq \frac{9C_0(1+\sigma \sqrt{d})}{\sqrt{np}},\label{eqn:12}\\
\opnorminf{\loo{U}\loo{H} - \loo{U^*}} &\leq 7 \br{ \opnorminf{UH - U^*}+\frac{1}{\sqrt{n}}}.   \label{eqn:15}
\end{align}
\end{lemma}

\begin{proof} See \proofref{lem:event-E-prime_new}.
\end{proof}

\subsection{Proof of Proposition \ref{prop:hat_H_new}}
We first state and prove a deterministic version of Proposition \ref{prop:hat_H_new}. The proof of Proposition \ref{prop:hat_H_new} follows from a simple probabilistic argument. %

\begin{lemma}\label{lem:hat_H}
Assume that (\ref{eqn:5}) holds. If $\frac{np}{(\sqrt{d} + \sigma d)^2} \geq 32 C_0^2$, we have
\begin{align}\label{eqn:11}
\min_{P\in\Pi_d}\fnorm{\hatH - PH^\top} \leq  \frac{{4C_0\br{\sqrt{d} + \sigma d}}}{\sqrt{np}}.
\end{align}
\end{lemma}
\begin{proof}
Note that
\begin{align*}
 \opnorm{U-U^*H^\top} &= \opnorm{UU^\top U - U^*U^{*\top}U} \leq  \opnorm{UU^\top  - U^*U^{*\top}}\opnorm{U}\leq \opnorm{UU^\top  - U^*U^{*\top}}.
\end{align*}
Since (\ref{eqn:5}) holds, we have
\begin{align*}
\fnorm{U - U^*H^\top} &\leq \sqrt{d}  \opnorm{U-U^*H^\top} \leq   \sqrt{d}\frac{8C_0\br{1+\sigma\sqrt{d}}}{7\sqrt{np}}.
\end{align*}

Note that $\sqrt{n}U^*$ has only $k$ unique rows $e_1^\top,\ldots,e_d^\top$ and each is of size $n$. Then $\sqrt{n}U^* H^\top$ also has $k$ unique rows $e_1^\top H^\top,\ldots,e_d^\top H^\top$ and each is of size $n$. From (\ref{eqn:5}), there exists a matrix $O\in \mathcal{O}_d$ such that $\opnorm{H-O}\leq \frac{8C_0\br{1+\sigma\sqrt{d}}}{7\sqrt{np}} \leq \frac{1}{7}$. Then
\begin{align*}
 \min_{a,b\in[d]} \norm{(e_a - e_b)^\top H^\top} &\geq \min_{a,b\in[d]} \norm{(e_a - e_b)^\top \br{O^\top+  (H-O)^\top}}\\
&\geq   \min_{a,b\in[d]} \norm{(e_a - e_b)^\top O^\top } - \opnorm{H-O}\\
&\geq \sqrt{2} -\frac{1}{7}\\
&>  \frac{6}{5}.
\end{align*}
That is, the minimum distance among the unique rows of $U^*H^\top$, is at least $\frac{6}{5\sqrt{n}}$.

Let $\hat z\in[d]^{nd}$ be the minimizer of (\ref{eqn:alg1}) with $\{\hat \mu_1,\ldots,\hat \mu_d\}$. Denote $\hat \Theta \define (\hat \mu_{\hat z_1}^\top, \ldots, \hat \mu_{\hat z_{nd}}^\top)^\top\in\mathr^{nd\times d}$. According to (\ref{eqn:alg1}), we have $\fnorm{\hat \Theta - U}\leq \fnorm{U^*H^\top - U}$. Hence,
\begin{align*}
\fnorm{\hat \Theta - U^*H^\top} &\leq \fnorm{\hat \Theta -U }  + \fnorm{U^*H^\top - U}\leq 2 \fnorm{U - U^*H^\top }\leq   \sqrt{d}\frac{16C_0\br{1+\sigma\sqrt{d}}}{7\sqrt{np}}.\numberthis \label{eqn:theta-err}
\end{align*}
Define $z^*\in[d]^{nd}$ such that the $i$th row of $U^* H^\top$ is equal to $e_{z^*_i}^\top H^T/\sqrt{n}$ for each $i\in[n]$.
Define the set $S$ as
\begin{align*}
S\define \cbr{i\in[nd]: \norm{\hat \mu_{\hat z_i} - e_{z^*_i}^\top H^T/\sqrt{n}} > \frac{3}{5\sqrt{n}}}.
\end{align*}
Then we have
\begin{align*}
\abs{S} \leq \frac{\fnorm{\hat \Theta - U^*H^\top}^2}{(\frac{3}{5\sqrt{n}})^2}\leq \frac{\br{4C_0\br{\sqrt{d} + \sigma d}}^2}{p}.
\end{align*}
Under the assumption $\frac{np}{(\sqrt{d} + \sigma d)^2} \geq 32 C_0^2,$
we have $\abs{S}\leq n/2$. Then
by the same argument as in the proof of Proposition 3.1 of \cite{zhang2022leave},  there exists a bijection $\phi:[d]\rightarrow[d]$ such that $\hat z_i = \phi(z^*_i)$ for all $i\notin S$. Hence, for each $a\in[d]$, we have
\begin{align*}
\norm{\hat \mu_{\phi(a)} - e_a^\top H^\top/\sqrt{n} }^2 & = \frac{\sum_{i\in [nd]:\hat z_i = \phi(a),z^*_i =a} \norm{\hat \mu_{\hat z_i} - e_{z^*_i}^\top H^\top }^2}{\abs{i\in[nd]:\hat z_i = \phi(a),z^*_i =a}} \leq   \frac{\sum_{i\in[nd]:\hat z_i = \phi(a)} \norm{\hat \mu_{\hat z_i} - e_{z^*_i}^\top H^\top }^2}{n-\abs{S}}\\
&\leq    \frac{\sum_{i\in[nd]:\hat z_i = \phi(a)} \norm{\hat \mu_{\hat z_i} - e_{z^*_i}^\top H^\top }^2}{n/2}.
\end{align*}
Hence,
\begin{align*}
\sum_{a\in[d]}\norm{\hat \mu_{\phi(a)} - e_a^\top H^\top/\sqrt{n}}^2 \leq \frac{\sum_{i\in[nd]} \norm{\hat \mu_{\hat z_i} - e_{z^*_i}^\top H^\top }^2 }{n/2} = \frac{\fnorm{\hat \Theta -U^* H^\top}^2}{n/2}  \leq \frac{\br{4C_0\br{\sqrt{d} + \sigma d}}^2}{n^2p}.
\end{align*}
That is, there exists a permutation matrix $P$ such that
\begin{align*}
\fnorm{\hatH - PH^\top}^2 \leq \frac{\br{4C_0\br{\sqrt{d} + \sigma d}}^2}{np},
\end{align*}
since rows of $M$  are $\sqrt{n}\hat \mu_1,\ldots,\sqrt{n}\hat \mu_d $.
\end{proof}

\begin{proof}[\textbf{Proof of Proposition \ref{prop:hat_H_new}}]
According to Lemma \ref{lem:event-E-prime_new}, there exist constants $C,C_0>0$ such that if $\frac{np}{\log n}>C$ and $\frac{np}{(1+\sigma \sqrt{d})^2 }\geq 64C_0^2$, we have   (\ref{eqn:5}) hold with probability at least $1-n^{-10}$. The proof is complete by Lemma \ref{lem:hat_H}.
\end{proof}

The following lemma establishes the accuracy of the anchor constructed using approximate clustering in Algorithm \ref{alg:approx-algo}. It serves as  a counterpart of Lemma \ref{lem:hat_H} and will be used in the proof of Theorem \ref{thm:main_approx}. Its proof is nearly identical to that of Lemma \ref{lem:hat_H} with minor modifications.
\begin{lemma} \label{lem:check_M}
Assume that (\ref{eqn:5}) holds. If $\frac{np}{(\sqrt{d} + \sigma d)^2} \geq 32 (1+\epsilon) C_0^2$, we have
\begin{align}\label{eqn:check_M}
\min_{P\in\Pi_d}\fnorm{\check M - PH^\top} \leq  \frac{{4C_0\sqrt{1+\epsilon}\br{\sqrt{d} + \sigma d}}}{\sqrt{np}}.
\end{align}
\end{lemma}
\begin{proof}
  We follow a similar argument to that in the proof of Lemma \ref{lem:hat_H}. The difference starts from the derivation of (\ref{eqn:theta-err}). 
  Define $\check \Theta = \left( \check \mu_{\check z_1}^\top, \ldots, \check\mu_{\tilde z_{nd}}^\top  \right)^\top $ to be the assigned centers given by the approximate clustering algorithm. Then, by (\ref{eqn:approx-cls-err}), we have $\fnorm{\check \Theta - U}^2\leq (1+\epsilon)\fnorm{U^*H^\top - U}^2$. Hence,
  \begin{align*}
    \fnorm{\check \Theta - U^*H^\top} &\leq \sqrt{1+\epsilon}\fnorm{\hat \Theta -U }  + \fnorm{U^*H^\top - U}\leq 2\sqrt{1+\epsilon} \fnorm{U - U^*H^\top }\\
    &\leq  \sqrt{1+\epsilon}\frac{16C_0\sqrt d\br{1+\sigma\sqrt{d}}}{7\sqrt{np}}\,.
  \end{align*}
The rest of the proof proceeds similarly, with a minor modification of the condition $\frac{np}{\br{\sqrt d + \sigma d}^2}\geq 32C_0^2$ to  $\frac{np}{\br{\sqrt d + \sigma d}^2}\geq 32C_0^2(1+\epsilon)$.
\end{proof}

\subsection{Proof of Theorem \ref{thm:l_infty_new}}

We first give a deterministic upper bound for $\opnormt{U_jH - U^*_j}$, using the decomposition (\ref{eqn:decomposition3}).
\begin{lemma}\label{lem:simplification}
Assume (\ref{eqn:2})-(\ref{eqn:12}) hold. Under the assumption $\frac{np}{(1+\sigma \sqrt{d})^2} \geq 22^2 C_0^2$, for each $j\in[n]$, we have
\begin{align}\label{eqn:13}
\opnorm{U_jH - U^*_j } &\leq  \frac{22C_0}{\sqrt{n}}\br{\frac{1+\sigma \sqrt{d}}{\sqrt{np}}}+2\opnorm{F_{j1}} + 2\opnorm{F_{j2}} + 2\opnorm{G_{j1}} + 2\opnorm{G_{j2}},
\end{align}
and
\begin{align}\label{eqn:Bj-bound}
\opnorm{B_j} &\leq  22C_0 \br{\frac{1}{\sqrt{n}}  +\opnorm{F_{j1}} + \opnorm{F_{j2}} + \opnorm{G_{j1}} + \opnorm{G_{j2}}} \br{\frac{1+\sigma \sqrt{d}}{\sqrt{np}}}.
\end{align}
\end{lemma}
\begin{proof}
Consider any $j\in[n]$. We have
\begin{align*}
\opnorm{B_j} &\leq  \opnorm{U_j }   \opnorm{H\Lambda - \Lambda H} \opnorm{\Lambda^{-1}}+ \opnorm{X_j}\opnorm{UU^\top - U^{(j)}U^{(j)\top}} \opnorm{\Lambda^{-1}}  \\
&\leq  {\frac{4}{3}\br{\opnorm{U_jH-U^*_j}+\frac{1}{\sqrt{n}}}} \br{2C_0\br{1+\sigma\sqrt{d}}\sqrt{np}} \frac{8}{7np} \\
&\quad +  \br{p\sqrt{n} + C_0\br{1+\sigma\sqrt{d}}\sqrt{np}} \br{6\br{\opnorm{U_jH -U^*_j}+\frac{1}{\sqrt{n}}}} \frac{8}{7np} \\
&\leq \frac{8}{7np}\br{\opnorm{U_jH-U^*_j}+\frac{1}{\sqrt{n}}} \br{6p\sqrt{n} + \frac{26}{3}C_0\br{1+\sigma \sqrt{d}}\sqrt{np}}\\
&\leq 11C_0 \br{\opnorm{U_jH-U^*_j}+\frac{1}{\sqrt{n}}} \br{\frac{1+\sigma \sqrt{d}}{\sqrt{np}}} ,
\end{align*}
where the second inequality is by (\ref{eqn:2})-(\ref{eqn:12}) and the last inequality is by $C_0 > 7$. From (\ref{eqn:decomposition3}), we also have
\begin{align*}
\opnorm{U_jH - U^*_j } &\leq \opnorm{B_j}+\opnorm{F_{j1}} + \opnorm{F_{j2}} + \opnorm{G_{j1}} + \opnorm{G_{j2}}.
\end{align*}
Plug the upper bound (\ref{eqn:Bj-bound}) of $\opnorm{B_j}$ into the above display. After rearrangement, we have
\begin{align*}
&\br{1- 11C_0 \br{\frac{1+\sigma \sqrt{d}}{\sqrt{np}}} }\opnorm{U_jH - U^*_j }\\
&\leq  \frac{11C_0}{\sqrt{n}}\br{\frac{1+\sigma \sqrt{d}}{\sqrt{np}}}+\opnorm{F_{j1}} + \opnorm{F_{j2}} + \opnorm{G_{j1}} + \opnorm{G_{j2}}.
\end{align*}
Under the assumption $\frac{np}{(1+\sigma \sqrt{d})^2} \geq 22^2 C_0^2$, we have
\begin{align*}
\opnorm{U_jH - U^*_j } &\leq \frac{22C_0}{\sqrt{n}}\br{\frac{1+\sigma \sqrt{d}}{\sqrt{np}}}+2\opnorm{F_{j1}} + 2\opnorm{F_{j2}} + 2\opnorm{G_{j1}} + 2\opnorm{G_{j2}}.
\end{align*}
Plugging it into the upper bound of $\opnorm{B_j}$, we have
\begin{align*}
&\opnorm{B_j}\\
&\leq 11C_0 \br{\br{\frac{22C_0}{\sqrt n}\br{\frac{1+\sigma \sqrt{d}}{\sqrt{np}}}+2\opnorm{F_{j1}} + 2\opnorm{F_{j2}} + 2\opnorm{G_{j1}} + 2\opnorm{G_{j2}}}+\frac{1}{\sqrt{n}}} \\
&\cdot \br{\frac{1+\sigma \sqrt{d}}{\sqrt{np}}} \\
&\leq  22C_0 \br{\frac{1}{\sqrt{n}}  +\opnorm{F_{j1}} + \opnorm{F_{j2}} + \opnorm{G_{j1}} + \opnorm{G_{j2}}} \br{\frac{1+\sigma \sqrt{d}}{\sqrt{np}}}.
\end{align*}
\end{proof}

\paragraph{Useful short-hand notations.} In our proofs, some key quantities are repeated and some are lexicographically cumbersome. In order to declutter the presentation, wherever convenient without sacrificing clarity, we will use the following shorthand notations. Define
\begin{align}
\Delta \define UH-U^*\in\mathr^{nd\times d},\quad  \Delta^{(j)} \define  U^{(j)}H^{(j)} - U^* \in\mathr^{nd\times d},\forall j\in[n], \label{eqn:Delta_def}
\end{align}
where $\loo{H} = \loo{U}^\top U^*$, such that block submatrices $\Delta_k = U_kH-U^*_k$ and $ \Delta^{(j)}_k = U^{(j)}_kH^{(j)} - U^*_k$ for each $j,k\in[n]$. 
We further introduce $\opnorminf{\cdot}$ norm such that
\begin{align}\label{eqn:Delta_norm_def}
\opnorminf{\Delta} &\define \max_{k\in[n]} \opnorm{\Delta_k},\quad\opnorminf{\Delta^{(j)}}\define \max_{k\in[n]} \opnorm{\Delta^{(j)}_k}.
\end{align}

\paragraph{Helper Tail Bound Inequalities and Proof Sketch.} As a preview, the terms $\opnorm{F_{j1}}$, $\opnorm{F_{j2}}$, $\opnorm{G_{j1}}$ and $\opnorm{G_{j2}}$ in (\ref{eqn:13}) and (\ref{eqn:Bj-bound}) can be bounded using the helper Lemma \ref{lem:Fj1}, Lemma \ref{lem:Fj2_denominator}, Lemma \ref{lem:Fj2} and Lemma \ref{lem:vershynin_new} below.  We defer the proof of these helper lemmas to a later section. Note that all of these inequalities require 
$\opnorminf{\Delta^{(j)}}$
which in turn can be bounded in terms of 
$\opnorminf{\Delta}$
using (\ref{eqn:15}). The key strategy is to construct an inequality in which 
$\opnorminf{\Delta}$
 appears on both the left and the right side of the inequality. After some manipulation and under appropriate assumptions, we reduce this to an inequality where the left side of the inequality consists of only 
 $\opnorminf{\Delta}$.

\begin{lemma}\label{lem:Fj1}
For any $j\in[n]$, both of the following tail bounds hold for any $t > 0$.
\begin{align*}
&\pbr{\opnorm{\sum_{k\neq j} A_{jk} \Delta_k^{(j)}}  \geq  p\sqrt{n} \opnorm{\Delta^{(j)}} + t\Bigg| \Delta^{(j)}} \\
&\leq 2d\ebr{- \frac{t^2/2}{p\sqrt{nd}\opnorm{\Delta^{(j)}}\opnorminf{\Delta^{(j)}} + \opnorminf{\Delta^{(j)}}t/3}}, \numberthis\label{eqn:Fj1-tail1}
\end{align*}
and
\begin{align*}
  &\pbr{\opnorm{\sum_{k\neq j} A_{jk} \Delta_k^{(j)}}  \geq  p\sqrt{n} \opnorm{\Delta^{(j)}} + t\Bigg| \Delta^{(j)}} \\
  &\leq 2d\ebr{- \frac{t^2/2}{np \opnorminf{\loo{\Delta}}^2 + \opnorminf{\Delta^{(j)}}t/3}}. \numberthis\label{eqn:Fj1-tail2}
  \end{align*}
\end{lemma}

\begin{proof} See \proofref{lem:Fj1}.
\end{proof}

\begin{lemma}\label{lem:Fj2_denominator}
For any $j\in[n]$, we have
\begin{align*}
&\pbr{\opnorm{ \sum_{k\neq j}A_{jk}(\Delta^{(j)}_k)^\top \Delta^{(j)}_k}   \geq  p \opnorm{\Delta^{(j)}}^2 + t \Bigg| \Delta^{(j)} }\\
& \leq 2d\ebr{ - \frac{t^2/2}{p \opnorminf{\Delta^{(j)}}^2  \opnorm{ \Delta^{(j)} }^2 + \opnorminf{\Delta^{(j)}}^2t/3}},
\end{align*}
for any $t\geq 0$.
\end{lemma}

\begin{proof} See \proofref{lem:Fj2_denominator}.
\end{proof}

\begin{lemma}\label{lem:Fj2}
There exists some constant $c>0$, such that for  any $t\geq 4\sqrt{d}\opnorm{ \sum_{k\neq j}A_{jk}(\Delta^{(j)}_k)^\top \Delta^{(j)}_k}^\frac{1}{2}$ and for any $j\in[n]$, we have
\begin{align*}
\pbr{\opnorm{\sum_{k\neq j}A_{jk}W_{jk}\Delta^{(j)}_k}\geq t\Bigg| \{A_{jk}\}_{k\neq j},\Delta^{(j)}}\leq 2\ebr{- \frac{ct^2}{\opnorm{ \sum_{k\neq j}A_{jk}(\Delta^{(j)}_k)^\top \Delta^{(j)}_k}}}.
\end{align*}
\end{lemma}

\begin{proof} See \proofref{lem:Fj2}.
\end{proof}

\begin{lemma}\label{lem:vershynin_new} [Corollary 7.3.3 of \cite{vershynin2018high}] Let $W$ be a $d\times d$ matrix with independent $\calN(0, 1)$ entries. Then there exists some constant $c>0$ such that  for every $t \geq 0$, we have
$$ \pbr{ \opnorm{W} \geq 2\sqrt{d} + t  } \leq 2\,\ebr{-ct^2} \,. $$
\end{lemma}

We are now ready to prove Theorem \ref{thm:l_infty_new}.

\begin{proof}[\textbf{Proof of Theorem \ref{thm:l_infty_new}}]
From Lemma \ref{lem:event-E-prime_new}, there exist constants $C_0', C_0 > 0$ such that if $\frac{np}{\log n} > C_0'$, then the event $\calE$ holds with probability at least $1-n^{-10}$. Assume $\calE$ holds and $\frac{np}{(\sqrt{d} + \sigma d)^2} \geq 22^2 C_0^2$. Then according to Lemma \ref{lem:event-E-prime_new}, we have (\ref{eqn:2})-(\ref{eqn:15}) hold. As a consequence, by Lemma \ref{lem:simplification}, we have (\ref{eqn:13}) hold as well. 

Consider any $j\in[n]$. Note that (\ref{eqn:13}) can be written as
\begin{align}
\opnorm{\Delta_j} &\leq  \frac{22C_0}{\sqrt{n}}\br{\frac{1+\sigma \sqrt{d}}{\sqrt{np}}}+2\opnorm{F_{j1}} + 2\opnorm{F_{j2}} + 2\opnorm{G_{j1}} + 2\opnorm{G_{j2}}.\label{eqn:17}
\end{align}
Hence, to upper bound $\opnorm{\Delta_j}$, 
we need to study $\opnorm{F_{j1}}$, $\opnorm{F_{j2}}$, $\opnorm{G_{j1}}$, and $\opnorm{G_{j2}}$.  

~\\
\textbf{Bounding $\opnorm{F_{j1}}$.} By (\ref{eqn:Fj1-tail2}) of Lemma \ref{lem:Fj1}, the following holds with probability at least $1-2dn^{-10}$
\begin{align*}
\opnorm{\sum_{k\neq j} A_{jk} \Delta_k^{(j)}}
&\leq   p\sqrt{n}   \opnorm{\Delta^{(j)}}  + \sqrt{40np\opnorminf{\Delta^{(j)}}^2 \log n}  + \frac{40}{3} \log n\opnorminf{\Delta^{(j)}}\\
&\leq   p\sqrt{n}   \opnorm{\Delta^{(j)}}  + \br{ \sqrt{40np \log n} + \frac{40}{3} \log n} \opnorminf{\Delta^{(j)}}.
\end{align*}
Then with (\ref{eqn:3}), we have
\begin{align*}
\opnorm{F_{j1}} &\leq \opnorm{\sum_{k\neq j} A_{jk} \Delta_k^{(j)}} \opnorm{\Lambda^{-1}}\\
 &\leq     \frac{8}{7np} \br{ p\sqrt{n}   \opnorm{\Delta^{(j)}}  + \br{ \sqrt{40np\log n} + \frac{40}{3} \log n} \opnorminf{\Delta^{(j)}}}\\
 &\leq   \frac{C_1}{\sqrt{n}} \opnorm{\Delta^{(j)}}  +C_1 \sqrt{\frac{\log n}{np}} \opnorminf{\Delta^{(j)}},\numberthis \label{eqn:18}
\end{align*}
for some constant $C_1>0$. 

~\\
\textbf{Bounding $\opnorm{F_{j2}}$.} By Lemma \ref{lem:Fj2}, we have
\begin{align*}
\opnorm{\sum_{k\neq j}A_{jk}W_{jk}\Delta^{(j)}_k} &\leq C_2 \sqrt{\log n} \opnorm{ \sum_{k\neq j}A_{jk}(\Delta^{(j)}_k)^\top \Delta^{(j)}_k}^\frac{1}{2}
\end{align*}
holds with probability at least $1-n^{-10}$ for some constant $C_2>0$. In addition, by Lemma \ref{lem:Fj2_denominator}, we have
\begin{align*}
\opnorm{ \sum_{k\neq j}A_{jk}(\Delta^{(j)}_k)^\top \Delta^{(j)}_k}  &\leq   p \opnorm{\Delta^{(j)}}^2 +\sqrt{40p\opnorm{\Delta^{(j)}}^2\opnorminf{\Delta^{(j)}}^2 \log n}  + \frac{40}{3} \log n\opnorminf{\Delta^{(j)}}^2\\
 &\leq   p \opnorm{\Delta^{(j)}}^2 +\sqrt{40np\opnorminf{\Delta^{(j)}}^4 \log n}  + \frac{40}{3} \log n\opnorminf{\Delta^{(j)}}^2\\
&\leq p \opnorm{\Delta^{(j)}}^2  +\br{ \sqrt{40np \log n} + \frac{40}{3} \log n} \opnorminf{\Delta^{(j)}}^2
\end{align*}
holds with probability at least $1-2dn^{-10}$.  Then with (\ref{eqn:3}), we have
\begin{align*}
\opnorm{F_{j2}} &\leq \sigma \opnorm{\sum_{k\neq j}A_{jk}W_{jk}\Delta^{(j)}_k}  \opnorm{\Lambda^{-1}}\\
&\leq  \frac{8\sigma}{7np}  C_2 \sqrt{\log n} \sqrt{p \opnorm{\Delta^{(j)}}^2  +\br{ \sqrt{40np \log n} + \frac{40}{3} \log n} \opnorminf{\Delta^{(j)}}^2}\\
&\leq  \frac{8\sigma}{7np}  C_2 \sqrt{\log n}   \br{\sqrt{p}  \opnorm{\Delta^{(j)}} + \sqrt{ \sqrt{40np \log n} + \frac{40}{3} \log n} \opnorminf{\Delta^{(j)}}}\\
&\leq  \frac{C_3\sigma}{\sqrt{n}}\sqrt{\frac{\log n}{np}} \opnorm{\Delta^{(j)}} + C_3\sigma \br{\frac{\log n}{np}}^\frac{3}{4}  \opnorminf{\Delta^{(j)}}, \numberthis \label{eqn:19}
\end{align*}
for some constant $C_3>0$. 

~\\
\textbf{For $\norm{G_{j1}}$}, its upper bound is given in (\ref{eqn:16}). 

~\\
\textbf{Bounding $\opnorm{G_{j2}}$.}
By Lemma \ref{lem:vershynin_new}, using the fact that $\{A_{jk}\}_{k\neq j}$ and $\{W_{jk}\}_{k\neq j}$ are independent, we have
\begin{align*}
\opnorm{\sum_{k\neq j} A_{jk} W_{jk}} &\leq \sqrt{\sum_{k\neq j} A_{jk}} \br{2\sqrt{d} + C_4\sqrt{\log n}},
\end{align*}
with probability at least $1-n^{-10}$,
for some constant $C_4>0$. Then
\begin{align*}
\opnorm{G_{j2}} &\leq \frac{\sigma}{\sqrt{n}} \opnorm{\sum_{k\neq j} A_{jk} W_{jk}} \opnorm{\Lambda^{-1}}\\
&\leq  \frac{\sigma}{\sqrt{n}}  \sqrt{\sum_{k\neq j} A_{jk}} \br{2\sqrt{d} + C_4\sqrt{\log n}}\opnorm{\Lambda^{-1}}\\
&\leq C_5\frac{\sigma}{\sqrt{n}}\br{\sqrt{\frac{d}{np}} + \sqrt{\frac{\log n }{np}}}, \numberthis \label{eqn:20}
\end{align*}
for some constant $C_5>0$, where the last inequality is by $\calE$ and (\ref{eqn:3}). 

~\\
\textbf{Putting things together.} Plugging (\ref{eqn:18}), (\ref{eqn:19}), (\ref{eqn:16}), and (\ref{eqn:20}) into (\ref{eqn:17}), we have
\begin{align*}
\opnorm{\Delta_j} &\leq  \frac{22C_0}{\sqrt{n}}\br{\frac{1+\sigma \sqrt{d}}{\sqrt{np}}}+2\br{ \frac{C_1}{\sqrt{n}} \opnorm{\Delta^{(j)}}  +C_1 \sqrt{\frac{\log n}{np}} \opnorminf{\Delta^{(j)}}} \\
&\quad + 2\br{ \frac{C_3\sigma}{\sqrt{n}}\sqrt{\frac{\log n}{np}} \opnorm{\Delta^{(j)}} + C_3\sigma \br{\frac{\log n}{np}}^\frac{3}{4}  \opnorminf{\Delta^{(j)}}} + \frac{4C_0}{\sqrt{n}}\br{\sqrt{\frac{\log n}{np}} + \frac{\sigma \sqrt{d}}{\sqrt{np}}} \\
&\quad + 2C_5\frac{\sigma}{\sqrt{n}}\br{\sqrt{\frac{d}{np}} + \sqrt{\frac{\log n }{np}}}\\
&\leq \frac{2}{\sqrt{n}} \br{C_1 + C_3\sigma\sqrt{\frac{\log n}{np}}} \opnorm{\Delta^{(j)}}  +2\br{C_1\sqrt{\frac{\log n}{np}} + C_3\sigma \br{\frac{\log n}{np}}^\frac{3}{4} }  \opnorminf{\Delta^{(j)}} \\
&\quad + \frac{26C_0 + 2C_5}{\sqrt{n}} \frac{\sqrt{\log n} + \sigma \sqrt{d} + \sigma\sqrt{\log n}}{\sqrt{np}}.
\end{align*}
By (\ref{eqn:15}) and (\ref{eqn:12}), we can replace the terms $\opnorm{\Delta^{(j)}}$, $\opnorminf{\Delta^{(j)}}$ with their respective upper bounds and have
\begin{align*}
\opnorm{\Delta_j} &\leq   \frac{2}{\sqrt{n}} \br{C_1 + C_3\sigma\sqrt{\frac{\log n}{np}}} \frac{9C_0(1+\sigma \sqrt{d})}{\sqrt{np}} \\
& +14\br{C_1\sqrt{\frac{\log n}{np}} + C_3\sigma \br{\frac{\log n}{np}}^\frac{3}{4} }   \br{ \opnorminf{\Delta}+\frac{1}{\sqrt{n}}} \\
&+ \frac{26C_0 + 2C_5}{\sqrt{n}} \frac{\sqrt{\log n} + \sigma \sqrt{d} + \sigma\sqrt{\log n}}{\sqrt{np}}.
\end{align*}
By a union bound, with  probability at least $1-6dn^{-9}$, the above inequality holds for all $j\in[n]$. Then
\begin{align*}
\opnorminf{\Delta} &\leq   \frac{2}{\sqrt{n}} \br{C_1 + C_3\sigma\sqrt{\frac{\log n}{np}}} \frac{9C_0(1+\sigma \sqrt{d})}{\sqrt{np}}\\
&+14\br{C_1\sqrt{\frac{\log n}{np}} + C_3\sigma \br{\frac{\log n}{np}}^\frac{3}{4} }   \br{ \opnorminf{\Delta}+\frac{1}{\sqrt{n}}} \\
& + \frac{26C_0 + 2C_5}{\sqrt{n}} \frac{\sqrt{\log n} + \sigma \sqrt{d} + \sigma\sqrt{\log n}}{\sqrt{np}}.
\end{align*}
After a rearrangement so that $\opnorminf{\Delta}$ only appears on the left side, we have
\begin{align*}
&\br{1- 14\br{C_1\sqrt{\frac{\log n}{np}} + C_3\sigma \br{\frac{\log n}{np}}^\frac{3}{4} }} \opnorminf{\Delta} \\
&\leq \frac{2}{\sqrt{n}} \br{C_1 + C_3\sigma\sqrt{\frac{\log n}{np}}} \frac{9C_0(1+\sigma \sqrt{d})}{\sqrt{np}} \\
&+ \frac{14}{\sqrt{n}}\br{C_1\sqrt{\frac{\log n}{np}} + C_3\sigma \br{\frac{\log n}{np}}^\frac{3}{4} }\\
&+ \frac{26C_0 + 2C_5}{\sqrt{n}} \frac{\sqrt{\log n} + \sigma \sqrt{d} + \sigma\sqrt{\log n}}{\sqrt{np}}.
\end{align*}
Then, there exists a constant $C>0$ such that if $\frac{np}{(\sigma^{4/3}  \vee 1)\log n} >C$, we have 
$$14\br{C_1\sqrt{\frac{\log n}{np}} + C_3\sigma \br{\frac{\log n}{np}}^\frac{3}{4} }\leq \frac{1}{2} .$$ Consequently, we have
\begin{align*}
\opnorminf{\Delta} &\leq \frac{4}{\sqrt{n}} \br{C_1 + C_3\sigma\sqrt{\frac{\log n}{np}}} \frac{9C_0(1+\sigma \sqrt{d})}{\sqrt{np}} +  \frac{28}{\sqrt{n}}\br{C_1\sqrt{\frac{\log n}{np}} + C_3 \sigma \br{\frac{\log n}{np}}^\frac{3}{4} }\\
&\quad + \frac{52C_0 + 4C_5}{\sqrt{n}} \frac{\sqrt{\log n} + \sigma \sqrt{d} + \sigma\sqrt{\log n}}{\sqrt{np}}\\
&\leq \frac{C_6}{\sqrt{n}}\br{\sqrt{\frac{\log n }{np}} + \frac{\sigma \sqrt{d}}{\sqrt{np}} + \frac{\sigma \sqrt{\log n}}{\sqrt{np}}},
\end{align*}
for some constant $C_6>0$.
\end{proof}

\subsection{Proofs of Theorem \ref{thm:main_new} and Theorem \ref{thm:main_approx}}

For the proof of Theorem \ref{thm:main_new}, we start with a proof sketch in Section \ref{subsubsec:1}. This sketch illustrates that the proof can be divided into four distinct steps. Subsequent detailed discussions of each step are elaborated upon in four separate sections, from Section \ref{subsubsec:2} to Section \ref{subsubsec:5}, with each section dedicated to one specific step of the proof. Following these, Section \ref{subsubsec:6} presents the proof of Theorem \ref{thm:main_approx}.

\subsubsection{Proof Sketch}\label{subsubsec:1}
In the previous section, the proof of Theorem 2 boils down to obtaining the operator norm bound on the terms $B_j, F_{j1}, F_{j2}, G_{j1}, G_{j2}$. In the proof of Theorem \ref{thm:main_new}, we go further by showing an \emph{exponential tail bound} on the operator norm of these terms. 

The proof of Theorem \ref{thm:main_new} is quite involved but can be broken down into the following key steps:

~\\
\emph{Step 1.} Let $\mathf$ (to be defined in (\ref{eqn:eventF})) be a high-probability event such that the expected value of the normalized Hamming loss $\E \ell(\hat Z,Z^*)\approx \E \ell(\hat Z,Z^*)\indic{\mathf}$. We first show an upper bound on $\E \ell(\hat Z,Z^*)\indic{\mathf}$ as a sum of $n$ terms, one for each block of $U$:
\begin{align}\label{eqn:new_4}
\E \ell(\hat Z,Z^*)\indic{\mathf} \leq \frac{1}{n}\sum_{j\in[n]} \sum_{R\in\Pi_d: R\neq I_d} \E \indic{\fnorm{  U_jM^\top - U^*_j \hat P^\top} \geq \fnorm{U_j M^\top - \frac{R \hat P^\top}{\sqrt{n}} }}\indic{\mathf},
\end{align}
where we recall the definition $\hat P \define \argmin_{P\in \Pi_d} \fnorm{M - PH^\top}$.
To understand it, note that each $\hat Z_j$ is obtained by finding a permutation matrix closest to $U_jM^\top$ as in (\ref{eqn:alg3}). If $\sqrt{n}U^*_j \hat P^\top$ is the closest permutation matrix to $U_jM^\top$, which can equivalently be stated as
\begin{align}\label{eqn:new_3}
\fnorm{  U_jM^\top - U^*_j \hat P^\top} < \fnorm{U_j M^\top - \frac{R \hat P^\top}{\sqrt{n}} },\forall R\in\Pi_d\text{ s.t. } R\neq I_d,
\end{align}
then according to (\ref{eqn:U_star_simplified}), $\hat Z_j = \sqrt{n}U^*_j \hat P^\top = \hat P^\top $. If the event (\ref{eqn:new_3}) holds for all $j\in[n]$, then all $\hat Z_j$ are $ \hat P^\top$, and consequently $\hat Z$ is equal to $Z^*$ up to a global permutation $ \hat P^\top$. Hence, the distance between $\hat Z$ and  $Z^*$ essentially depends on the number of $j\in[n]$ for which (\ref{eqn:new_3}) is not satisfied, which can be further upper bounded and leads to (\ref{eqn:new_4}).

~\\
\emph{Step 2.} From (\ref{eqn:new_3}), upper bounding $\E \ell(\hat Z,Z^*)\indic{\mathf}$ is about upper bounding each individual term on its right-hand side. It turn out each  error term  can be further decomposed into four tail probabilities. Consider any $j\in[n]$ and any $R\in\Pi_d$ such that $R\neq I_d$, we are able to show
\begin{align*}
& \E \indic{\fnorm{  U_jM^\top - U^*_j \hat P^\top} \geq \fnorm{U_j M^\top - \frac{R \hat P^\top}{\sqrt{n}} }}\indic{\mathf} \\
&\leq \text{ a tail probability of }\opnorm{F_{j1}} \quad + \quad \text{ a tail probability of }\opnorm{F_{j2}} \\
 &\quad + \text{ a tail probability of }\opnorm{G_{j2}} \quad + \quad \text{ a tail probability of }\iprod{\sum_{k\neq j} A_{jk} W_{jk} }{R-I_d} . \numberthis \label{eqn:new_5}
\end{align*}
This step is summarized in Lemma \ref{lem:hamming-decomp} and explicit expressions of tail probabilities are shown in (\ref{eqn:hamming-decomp-ineq}). Among the four tail probabilities, the last one that involves $\sum_{k\neq j} A_{jk} W_{jk}$ will lead to the exponential error rate. There is a variable $\rho$ that trade-offs the contributions between the four terms which will be determined in the last step. 

~\\
\emph{Step 3.} In this step, we leverage Lemma \ref{lem:Fj1}-\ref{lem:vershynin_new} to obtain exponential bound on each of the four tail probabilities of (\ref{eqn:new_5}).  Despite $\opnorm{F_{j1}}$, $\opnorm{F_{j2}} $, and $\opnorm{G_{j2}}$ have been analyzed in the proof of Theorem \ref{thm:l_infty_new}, results established there can not directly applied here. In fact, to derive sharp tail bounds here, we need leverage the obtained blockwise maximum deviation in Theorem \ref{thm:l_infty_new}. 
The final tail bounds are given in (\ref{eqn:31}), (\ref{eqn:32}), (\ref{eqn:33}),  and (\ref{eqn:34}), respectively.

~\\
\emph{Step 4.}
In the last step, we combine the four upper bounds in Step 3 to obtain a tight upper bound on (\ref{eqn:new_5}), and consequently on $\E \ell(\hat Z,Z^*)\indic{\mathf}$. At this point, we will select the variable $\rho$ introduced in Step 2 to ensure the tail probability involving $\sum_{k\neq j} A_{jk} W_{jk}$ dominates the other three, and consequently
 the final error bound matches the desired form.

\subsubsection{Step 1. Decomposition of The Hamming Loss}  \label{subsubsec:2}

From Lemma \ref{lem:event-E-prime_new}, there exist constants $C_0', C_0 > 0$ such that if $\frac{np}{\log n} > C_0'$, then the event $\calE$ holds with probability at least $1-n^{-10}$. By Theorem \ref{thm:l_infty_new}, there exist some constants $C,C',C''>0$ such that if $\frac{np}{(\sigma^{4/3}  \vee 1)\log n} >C$ and $\frac{np}{(\sqrt{d} + \sigma d)^2}\geq C'$, we have (\ref{eqn:22}) holds with probability at least  $1-6dn^{-9}$. Denote $\mathf$ to be the event that both $\calE$ and (\ref{eqn:22}) hold. That is,
\begin{equation}\label{eqn:eventF}
\mathf \define \left\{ \calE, (\ref{eqn:22}) \text{ holds}  \right\}\,.
\end{equation}
By union bound, we have
\begin{align}\label{eqn:36}
\pbr{\mathf} \geq 1-7dn^{-9}.
\end{align}
Assume $\mathf$ holds and $\frac{np}{(\sqrt{d} + \sigma d)^2} \geq 22^2 C_0^2$. Then according to Lemma \ref{lem:event-E-prime_new}, we have (\ref{eqn:2})-(\ref{eqn:15}) hold. As a consequence, by Lemma \ref{lem:simplification}  and Lemma \ref{lem:hat_H}, we have (\ref{eqn:13})-(\ref{eqn:Bj-bound}) and (\ref{eqn:11}) all hold as well. Note that 
\begin{align}
\E\ell(\hat Z,Z^*) \leq \E\ell(\hat Z,Z^*) \indic{\mathf}  + \pbr{\mathf^\mathrm{c}} \leq \E\ell(\hat Z,Z^*) \indic{\mathf} + 7dn^{-9} \label{eqn:35}\,.
\end{align}
We will later show that the $7dn^{-9}$ is dominated by $ \E\ell(\hat Z,Z^*) \indic{\mathf} $. Define
\begin{align*}
    \hat P \define \argmin_{P\in \Pi_d} \, \fnorm{\hatH - P H^\top}.
\end{align*}
Since we let $Z^*_j=I_d$ for all $j\in[n]$, we have
\begin{align*}
&\E \ell(\hat Z,Z^*)  \indic{\mathf} \\
&\leq \E \br{ \frac{1}{n} \sum_{j\in[n]} \indic{\hat Z_j \neq Z^*_j \hat P^\top}  \indic{\mathf} } = \frac{1}{n}\sum_{j\in[n]} \E  \indic{\hat Z_j \neq Z^*_j \hat P^\top}  \indic{\mathf} \\
&\leq \frac{1}{n}\sum_{j\in[n]} \E  \indic{\exists R\in\Pi_d \text{ s.t. }R\neq I_d \text{ and }\fnorm{ \sqrt{n} U_jM^\top - R \hat P^\top }^2 \leq \fnorm{  \sqrt{n} U_jM^\top - Z^*_j \hat P^\top}^2}  \indic{\mathf} \\
& =  \frac{1}{n}\sum_{j\in[n]} \E  \indic{\exists R\in\Pi_d \text{ s.t. }R\neq I_d \text{ and }\fnorm{  U_jM^\top - \frac{R \hat P^\top}{\sqrt{n}} }^2 \leq \fnorm{  U_jM^\top - U^*_j \hat P^\top}^2}  \indic{\mathf}\\
&\leq  \frac{1}{n}\sum_{j\in[n]} \sum_{R\in\Pi_d: R\neq I_d} \E \indic{\fnorm{U_j M^\top - \frac{R \hat P^\top}{\sqrt{n}} }^2 \leq \fnorm{  U_jM^\top - U^*_j \hat P^\top}^2}  \indic{\mathf} \\
&= \frac{1}{n}\sum_{j\in[n]} \sum_{R\in\Pi_d: R\neq I_d} \E \indic{\fnorm{  U_jM^\top - U^*_j \hat P^\top} \geq \fnorm{U_j M^\top - \frac{R \hat P^\top}{\sqrt{n}} }}  \indic{\mathf}, \numberthis \label{eqn:26}
\end{align*}
where the second inequity is by the definition of $\hat Z_j$ in (\ref{eqn:alg3}) and the second equation is by that we let $Z^*_j=I_d$.

\subsubsection{Step 2. Decomposition of Each Individual Error} 
The reader can see that the key to obtaining our final error bound is to obtain a tight upper bound on (\ref{eqn:26}), which is all about analyzing each individual error $ \E \indic{\fnorm{  U_jM^\top - U^*_j \hat P^\top} \geq \fnorm{U_j M^\top - \frac{R \hat P^\top}{\sqrt{n}} }}  \indic{\mathf}$.
The next step is the following inequality which shows an upper bound on it consists in four terms. Note that we introduce a $\rho$ parameter which trades off the contribution among the three terms. $\rho$ will be carefully chosen in the third step. Recall the decomposition (\ref{eqn:decomposition3}).

\begin{lemma}\label{lem:hamming-decomp} Consider any $j\in[n]$ and any $R \in \Pi_d$ such that $R\neq I_d$. Define
\begin{align*}
h_R \define \fnorm{R - I_d}^2/2.
\end{align*} 
Suppose that (\ref{eqn:16}), (\ref{eqn:10}), (\ref{eqn:11}), (\ref{eqn:Bj-bound}), and $\frac{np}{(1+\sigma \sqrt{d})^2} \geq 22^2 C_0^2$ hold. Then there exists some $C_2>0$ such that for any $\rho > 0$,
the following inequality holds deterministically:
\begin{align*}
  &\E\indic{ \fnorm{U_j \hatH^\top - U_j^*\hat P^\top} \geq \fnorm{U_j \hatH^\top - \frac{R \hat P^\top}{\sqrt n} }   }  \indic{\mathf} \\
  &\leq \E\indic{4\sqrt{2dh_R} \opnorm{F_{j1}}  \geq \rho \frac{h_R}{\sqrt{n}}}\indic{\mathf}  \\
  &+ \E\indic{4\sqrt{2dh_R} \opnorm{F_{j2}}  \geq \rho \frac{h_R}{\sqrt{n}}}\indic{\mathf}\\
  &+ \E \indic{ 46C_0\br{\frac{1+\sigma \sqrt{d}}{\sqrt{np}}} \sqrt{2dh_R} \opnorm{G_{j2}} \geq  \rho \frac{h_R}{\sqrt{n}}}\indic{\mathf}   \\
  & + \E \indic{ \frac{\sigma}{\sqrt{n}}\frac{1}{np}  \iprod{\sum_{k\neq j} A_{jk} W_{jk} }{R-I_d} \geq \br{1- 3\rho - C_2\sqrt d\br{\frac{\sqrt{\log n}+ \sigma\sqrt{d}}{\sqrt{np}}}}     \frac{h_R}{\sqrt{n}} }\indic{\mathf}.  \numberthis \label{eqn:hamming-decomp-ineq}
  \end{align*}
\end{lemma}

\begin{proof} 
By the definition of $h_R$, we have $2\leq h_R\leq d$. We have
\begin{align*}
&\indic{ \fnorm{U_j \hatH^\top - U_j^*\hat P^\top} \geq \fnorm{U_j \hatH^\top - \frac{R \hat P^\top}{\sqrt n} }   }  \indic{\mathf} \\
&= \indic{\iprod{U_j \hatH^\top - U_j^*\hat P^\top}{ \frac{R \hat P^\top}{\sqrt n} - U_j^*\hat P^\top}\geq \frac{1}{2}\fnorm{ \frac{R \hat P^\top}{\sqrt n} - U_j^*\hat P^\top}^2 } \indic{\mathf} \\
&=\indic{\iprod{U_j \hatH^\top \hat P - U_j^*}{ \frac{R}{\sqrt n} - U_j^*}\geq \frac{1}{2}\fnorm{ \frac{R }{\sqrt n} - U_j^*}^2 } \indic{\mathf}  \\
& = \indic{\iprod{U_j \hatH^\top \hat P - U_j^*}{R-I_d} \geq  \frac{h_R}{\sqrt{n}}} \indic{\mathf}  \\
& = \indic{\iprod{U_j H - U_j^*}{R-I_d} + \iprod{U_j\br{M - \hat PH^\top}^\top \hat P}{ R-I_d} \geq  \frac{h_R}{\sqrt{n}}} \indic{\mathf}. \numberthis \label{eqn:9}
\end{align*}

We have to obtain a more refined decomposition of both inner products in (\ref{eqn:9}). For the second inner product of (\ref{eqn:9}) we have
\begin{align*}
&\iprod{U_j\br{M - \hat PH^\top}^\top \hat P}{ R-I_d} \\
&\leq  \fnorm{U_j \br{M -\hat P  H^\top}^\top}\fnorm{ R-I_d}\\
&\leq  \opnorm{U_j} \fnorm{M -\hat P  H^\top}\fnorm{ R-I_d}\\
&\leq   \frac{4}{3}\br{\opnorm{U_jH-U^*_j}+\frac{1}{\sqrt{n}}}\br{ \frac{{4C_0\br{\sqrt{d} + \sigma d}}}{\sqrt{np}}}\sqrt{2h_R}\\
&\leq  \frac{4}{3}\br{\opnorm{B_j}+\opnorm{F_{j1}}+\opnorm{F_{j2}}+\opnorm{G_{j1}}+ \opnorm{G_{j2}}+\frac{1}{\sqrt{n}}}\br{ \frac{{4C_0\br{\sqrt{d} + \sigma d}}}{\sqrt{np}}}\sqrt{2h_R},\numberthis \label{eqn:second-inner}
\end{align*}
where the third inequality is by (\ref{eqn:10}) and (\ref{eqn:11}) and the fourth inequality is by (\ref{eqn:decomposition3}).

For the first inner product of (\ref{eqn:9}), by (\ref{eqn:decomposition3}), we have
\begin{align*}
\iprod{U_j H - U_j^*}{R-I_d} &= \iprod{B_j+ F_{j1} + F_{j2} + G_{j1} }{R-I_d}    + \iprod{G_{j2}}{R-I_d}\\
&\leq \fnorm{B_j+ F_{j1} + F_{j2} + G_{j1}} \fnorm{R-I_d}  + \iprod{G_{j2}}{R-I_d}\\
&\leq \sqrt{d}\opnorm{B_j+ F_{j1} + F_{j2} + G_{j1}}  \sqrt{2h_R} + \iprod{G_{j2}}{R-I_d}\\
&\leq  \br{\opnorm{B_j} + \opnorm{F_{j1}} + \opnorm{F_{j2}} + \opnorm{G_{j1}} } \sqrt{2dh_R} + \iprod{G_{j2}}{R-I_d}.
\end{align*}
From (\ref{eqn:Bj-bound}), since $\frac{np}{(1+\sigma \sqrt{d})^2} \geq 22^2 C_0^2$, we have
\begin{align}\label{eqn:new_2}
\opnorm{B_j} &\leq  \frac{22C_0}{\sqrt{n}} \br{\frac{1+\sigma \sqrt{d}}{\sqrt{np}}} + \opnorm{F_{j1}} + \opnorm{F_{j2}} + \opnorm{G_{j1}}  + 22C_0\br{\frac{1+\sigma \sqrt{d}}{\sqrt{np}}}  \opnorm{G_{j2}}.
\end{align}
The reader can see that thanks to this upper bound, we only have to contend with establishing upper bounds on the operator norm of $F_{j1}, F_{j2}, G_{j1}, G_{j2}$. From the above three displays, we proceed with a series of inequalities where we use the established upper bound results obtained previously to obtain an upper bound that can be written in terms of $\opnorm{F_{j1}}, \opnorm{F_{j2}}, \opnorm{G_{j1}}$ and $\iprod{G_{j2}}{R- I_d}$.
\begin{align*}
&\iprod{U_j H - U_j^*}{R-I_d} + \iprod{U_j\br{M - \hat PH^\top}^\top \hat P}{ R-I_d} \\
&\leq \br{1+ \frac{16C_0(1+\sigma\sqrt{d})}{3\sqrt{np}} }\br{\opnorm{B_j} + \opnorm{F_{j1}} + \opnorm{F_{j2}} + \opnorm{G_{j1}} } \sqrt{2dh_R} \\
&\quad +\frac{16C_0(1+\sigma\sqrt{d})}{3\sqrt{np}} \sqrt{2dh_R}\opnorm{G_{j2}}+ \iprod{G_{j2}}{R-I_d} + \frac{16C_0(1+\sigma\sqrt{d})}{3\sqrt{np}} \frac{\sqrt{2dh_R}}{\sqrt{n}}\\
&\leq  2\br{1+ \frac{16C_0(1+\sigma\sqrt{d})}{3\sqrt{np}} }\br{ \opnorm{F_{j1}} + \opnorm{F_{j2}} + \opnorm{G_{j1}} } \sqrt{2dh_R} \\
&\quad +\br{\frac{41}{33} +\frac{16C_0(1+\sigma\sqrt{d})}{3\sqrt{np}}  }22C_0\br{\frac{1+\sigma \sqrt{d}}{\sqrt{np}}}  \opnorm{G_{j2}} \sqrt{2dh_R} + \iprod{G_{j2}}{R-I_d}\\
&\quad +\br{\frac{41}{33} +\frac{16C_0(1+\sigma\sqrt{d})}{3\sqrt{np}}  }22C_0\br{\frac{1+\sigma \sqrt{d}}{\sqrt{np}}} \frac{\sqrt{2dh_R} }{\sqrt{n}}.
\end{align*}
where the second inequality is by the upper bound on $\opnorm{B_{j}}$ from (\ref{eqn:new_2}). We then have
\begin{align*}
&\iprod{U_j H - U_j^*}{R-I_d} + \iprod{U_j\br{M - \hat PH^\top}^\top \hat P}{ R-I_d} \\
&\leq  2\br{1+ \frac{16C_0(1+\sigma\sqrt{d})}{3\sqrt{np}} }\br{ \opnorm{F_{j1}} + \opnorm{F_{j2}} } \sqrt{2dh_R} \\
&\quad +\br{\frac{41}{33} +\frac{16C_0(1+\sigma\sqrt{d})}{3\sqrt{np}}  }22C_0\br{\frac{1+\sigma \sqrt{d}}{\sqrt{np}}}  \opnorm{G_{j2}} \sqrt{2dh_R} + \iprod{G_{j2}}{R-I_d}\\
&\quad +\br{\frac{41}{33} +\frac{16C_0(1+\sigma\sqrt{d})}{3\sqrt{np}}  }22C_0\br{\frac{1+\sigma \sqrt{d}}{\sqrt{np}}} \frac{\sqrt{2dh_R} }{\sqrt{n}}\\
&\quad + 4\br{1+ \frac{16C_0(1+\sigma\sqrt{d})}{3\sqrt{np}} } C_0\br{\frac{\sqrt{\log n} + \sigma\sqrt{d}}{\sqrt{np}}} \frac{\sqrt{2dh_R}}{\sqrt{n}}\\
&\leq  2\br{1+ \frac{16C_0(1+\sigma\sqrt{d})}{3\sqrt{np}} }\br{ \opnorm{F_{j1}} + \opnorm{F_{j2}} } \sqrt{2dh_R} \\
&\quad +\br{\frac{41}{33} +\frac{16C_0(1+\sigma\sqrt{d})}{3\sqrt{np}}  }22C_0\br{\frac{1+\sigma \sqrt{d}}{\sqrt{np}}}  \opnorm{G_{j2}} \sqrt{2dh_R} + \iprod{G_{j2}}{R-I_d}\\
&\quad +\br{\frac{41}{33} +\frac{16C_0(1+\sigma\sqrt{d})}{3\sqrt{np}}  }26C_0\br{\frac{\sqrt{\log n}+\sigma \sqrt{d}}{\sqrt{np}}} \frac{\sqrt{2dh_R} }{\sqrt{n}}\,,
\end{align*}
where the first inequality is by the upper bound on $\opnorm{G_{j1}}$ in (\ref{eqn:16}) and the last inequality is obtained by combining the two terms containing $\frac{\sqrt{2dh_R}}{\sqrt n}$.

In addition to $\opnorm{F_{j1}}$, $\opnorm{F_{j2}}$, and $\opnorm{G_{j1}}$, the above display involves 
$ \iprod{G_{j2}}{R-I_d}$ which must be analyzed closely in order to achieve the optimal error bound. 
Note the following decomposition.
\begin{align*}
&\iprod{G_{j2}}{R-I_d} \\
& = \frac{\sigma}{\sqrt{n}} \iprod{\sum_{k\neq j} A_{jk} W_{jk} \Lambda^{-1}}{R-I_d} \\
&= \frac{\sigma}{\sqrt{n}} \iprod{\sum_{k\neq j} A_{jk} W_{jk} \frac{1}{np} }{R-I_d}  + \frac{\sigma}{\sqrt{n}} \iprod{\sum_{k\neq j} A_{jk} W_{jk} \br{\Lambda^{-1} -  \frac{1}{np}I_d} }{R-I_d} \\
&= \frac{\sigma}{\sqrt{n}} \frac{1}{np} \iprod{\sum_{k\neq j} A_{jk} W_{jk} }{R-I_d}  + \frac{\sigma}{\sqrt{n}} \iprod{\sum_{k\neq j} A_{jk} W_{jk} \Lambda^{-1}(npI_d-\Lambda) \frac{1}{np} }{R-I_d} \\
&= \frac{\sigma}{\sqrt{n}}  \frac{1}{np}\iprod{\sum_{k\neq j} A_{jk} W_{jk} }{R-I_d}  +  \frac{1}{np}\iprod{G_{j2}(npI_d-\Lambda) }{R-I_d} \\
&\leq  \frac{\sigma}{\sqrt{n}}  \frac{1}{np}\iprod{\sum_{k\neq j} A_{jk} W_{jk} }{R-I_d}  +  \frac{\sqrt{d}}{np}\opnorm{G_{j2} (npI_d -\Lambda)} \fnorm{R-I_d}\\
&\leq   \frac{\sigma}{\sqrt{n}}\frac{1}{np} \iprod{\sum_{k\neq j} A_{jk} W_{jk} }{R-I_d}  + \frac{\sqrt{d}}{np}\opnorm{G_{j2} }\max_{i\in[n]} \abs{\lambda_i(X)-np} \sqrt{2h_R}\\
&\leq   \frac{\sigma}{\sqrt{n}}\frac{1}{np}  \iprod{\sum_{k\neq j} A_{jk} W_{jk} }{R-I_d}  + \frac{C_0\sqrt{2dh_R}(1+\sigma\sqrt{d})}{\sqrt{np}} \opnorm{G_{j2} }.
\end{align*}
Hence, plugging the above inequality into the display for 
$$\iprod{U_j H - U_j^*}{R-I_d} + \iprod{U_j\br{M - \hat PH^\top}^\top \hat P}{ R-I_d}$$ 
gives
\begin{align*}
&\iprod{U_j H - U_j^*}{R-I_d} + \iprod{U_j\br{M - \hat PH^\top}^\top \hat P}{ R-I_d} \\
&\leq  2\br{1+ \frac{16C_0(1+\sigma\sqrt{d})}{3\sqrt{np}} }\br{ \opnorm{F_{j1}} + \opnorm{F_{j2}} } \sqrt{2dh_R} \\
&\quad +\br{\frac{41}{33} +\frac{16C_0(1+\sigma\sqrt{d})}{3\sqrt{np}}  }23C_0\br{\frac{1+\sigma \sqrt{d}}{\sqrt{np}}}  \opnorm{G_{j2}} \sqrt{2dh_R} \\
&\quad +\br{\frac{41}{33} +\frac{16C_0(1+\sigma\sqrt{d})}{3\sqrt{np}}  }26C_0\br{\frac{\sqrt{\log n}+\sigma \sqrt{d}}{\sqrt{np}}} \frac{\sqrt{2dh_R} }{\sqrt{n}}\\
&\quad +  \frac{\sigma}{\sqrt{n}}\frac{1}{np}  \iprod{\sum_{k\neq j} A_{jk} W_{jk} }{R-I_d}\\
&\leq 4\br{ \opnorm{F_{j1}} + \opnorm{F_{j2}} } \sqrt{2dh_R}  + 46C_0\br{\frac{1+\sigma \sqrt{d}}{\sqrt{np}}}  \opnorm{G_{j2}} \sqrt{2dh_R} \\
&\quad + 52C_0\br{\frac{\sqrt{\log n}+\sigma \sqrt{d}}{\sqrt{np}}} \frac{\sqrt{2dh_R} }{\sqrt{n}} +  \frac{\sigma}{\sqrt{n}}\frac{1}{np}  \iprod{\sum_{k\neq j} A_{jk} W_{jk} }{R-I_d}, \numberthis \label{eqn:25}
\end{align*}
where the last inequality holds when  $\frac{\sqrt{np}}{1+\sigma\sqrt{d}} >\frac{176}{25}C_0$.

As a result, by (\ref{eqn:25}), (\ref{eqn:9}) becomes
\begin{align*}
&\indic{ \fnorm{U_j \hatH^\top - U_j^*\hat P^\top} \geq \fnorm{U_j \hatH^\top - \frac{R \hat P^\top}{\sqrt n} } }  \indic{\mathf} \\
&\leq   \mathbb{I}\Bigg\{ 4\sqrt{2dh_R} \br{\opnorm{F_{j1}} + \opnorm{F_{j2}}  } +   46C_0\br{\frac{1+\sigma \sqrt{d}}{\sqrt{np}}} \sqrt{2dh_R} \opnorm{G_{j2}}\\
&\quad + \frac{\sigma}{\sqrt{n}}\frac{1}{np}  \iprod{\sum_{k\neq j} A_{jk} W_{jk} }{R-I_d} \geq \br{1- C_2 \sqrt d\br{\frac{\sqrt{\log n}+ \sigma\sqrt{d}}{\sqrt{np}}}}     \frac{h_R}{\sqrt{n}} \Bigg\} \indic{\mathf},
\end{align*}
for some constant $C_2>0$.  Let $\rho>0$ be some quantity whose value will be determined later. We have
\begin{align*}
&\indic{ \fnorm{U_j \hatH^\top - U_j^*\hat P^\top} \geq \fnorm{U_j \hatH^\top - \frac{R \hat P^\top}{\sqrt n} }   }  \indic{\mathf} \\
&\leq \indic{4\sqrt{2dh_R} \opnorm{F_{j1}}  \geq \rho \frac{h_R}{\sqrt{n}}}\indic{\mathf}  + \indic{4\sqrt{2dh_R} \opnorm{F_{j2}}  \geq \rho \frac{h_R}{\sqrt{n}}}\indic{\mathf}\\
&\quad +  \indic{ 46C_0\br{\frac{1+\sigma \sqrt{d}}{\sqrt{np}}} \sqrt{2dh_R} \opnorm{G_{j2}} \geq  \rho \frac{h_R}{\sqrt{n}}}\indic{\mathf}   \\
&\quad +  \indic{ \frac{\sigma}{\sqrt{n}}\frac{1}{np}  \iprod{\sum_{k\neq j} A_{jk} W_{jk} }{R-I_d} \geq \br{1- 3\rho - C_2\sqrt d\br{\frac{\sqrt{\log n}+ \sigma\sqrt{d}}{\sqrt{np}}}}     \frac{h_R}{\sqrt{n}} }\indic{\mathf}. 
\end{align*}
Taking expected values on both sides, we complete the proof.
\end{proof}

\subsubsection{Step 3. Bounding The Four Tail Probabilities}
In the following, we are going to establish upper bounds for each of the four tail probabilities in (\ref{eqn:hamming-decomp-ineq}). Recall the definitions of $\Delta,\Delta^{(j)}$ in (\ref{eqn:Delta_def}) and those of  $\opnorminf{\Delta},\opnorminf{\Delta^{(j)}}$ in (\ref{eqn:Delta_norm_def}).

\textbf{For the first term in (\ref{eqn:hamming-decomp-ineq})}, by (\ref{eqn:3}), we have
\begin{align*}
\E\indic{4\sqrt{2dh_R} \opnorm{F_{j1}}  \geq \rho \frac{h_R}{\sqrt{n}}}\indic{\mathf} &\leq  \E\indic{4\sqrt{2dh_R}\opnorm{\sum_{k\neq j} A_{jk} \Delta^{(j)}_k} \opnorm{\Lambda^{-1}}  \geq \rho \frac{h_R}{\sqrt{n}}}\indic{\mathf} \\
 &\leq \E\indic{\frac{32\sqrt{2dh_R} }{7np}\opnorm{\sum_{k\neq j} A_{jk} \Delta^{(j)}_k}  \geq \rho \frac{h_R}{\sqrt{n}}}\indic{\mathf}.
\end{align*}
Define an event
\begin{align*}
\mathf_{j}\define \mathbb{I} \Bigg\{ & \opnorm{\Delta^{(j)}} \leq \frac{9C_0(1+\sigma \sqrt{d})}{\sqrt{np}},\\
&\opnorminf{\Delta^{(j)}}\leq 7\br{\frac{ C''}{\sqrt{n}}\br{\sqrt{\frac{\log n }{np}} + \frac{\sigma \sqrt{d}}{\sqrt{np}} + \frac{\sigma \sqrt{\log n}}{\sqrt{np}}} + \frac{1}{\sqrt{n}}}
\Bigg\}.
\end{align*}
Note that the upper bound for $\opnorminf{\Delta^{(j)}}$ in $\mathf_j$ is a direct consequence of (\ref{eqn:15}) and (\ref{eqn:22}). Together with (\ref{eqn:12}), we have $\mathf\subset \mathf_j$. The event $\mathf_j$ has an equivalent expression. Define a set
\begin{align*}
\mathcal{Y}\define \Bigg\{ Y&=(Y_1^\top,Y_2^\top,\ldots,Y_n^\top)^\top\in\mathr^{nd\times d}: \opnorm{Y} \leq \frac{9C_0(1+\sigma \sqrt{d})}{\sqrt{np}},\\
&\max_{i\in[n]}\opnorm{Y_i}\leq 7\br{\frac{ C''}{\sqrt{n}}\br{\sqrt{\frac{\log n }{np}} + \frac{\sigma \sqrt{d}}{\sqrt{np}} + \frac{\sigma \sqrt{\log n}}{\sqrt{np}}} + \frac{1}{\sqrt{n}}}
\Bigg\}.
\end{align*}
Then $\mathf_j = \indic{\Delta^{(j)} \in \mathcal{Y}}$. Using that $\Delta^{(j)}$ is independent of $\{A_{jk}\}_{k\neq j}$,
we have
\begin{align*}
&\E\indic{4\sqrt{2dh_R} \opnorm{F_{j1}}  \geq \rho \frac{h_R}{\sqrt{n}}}\indic{\mathf} \\
&\leq \E\indic{\frac{32\sqrt{2dh_R} }{7np}\opnorm{\sum_{k\neq j} A_{jk} \Delta^{(j)}_k}  \geq \rho \frac{h_R}{\sqrt{n}}}\indic{\mathf_j}\\
& = \pbr{\frac{32\sqrt{2dh_R} }{7np}\opnorm{\sum_{k\neq j} A_{jk} \Delta^{(j)}_k}  \geq \rho \frac{h_R}{\sqrt{n}}, \mathf_j}\\
& = \E \br{  \pbr{\frac{32\sqrt{2dh_R} }{7np}\opnorm{\sum_{k\neq j} A_{jk} \Delta^{(j)}_k}  \geq \rho \frac{h_R}{\sqrt{n}}, \mathf_j \Bigg| \Delta^{(j)}}}\\
& = \E\br{  \pbr{\frac{32\sqrt{2dh_R} }{7np}\opnorm{\sum_{k\neq j} A_{jk} \Delta^{(j)}_k}  \geq \rho \frac{h_R}{\sqrt{n}} \Bigg| \Delta^{(j)}} \indic{ \Delta^{(j)} \in \mathf_j}}\\
&\leq  \sup_{\Delta^{(j)}\in\mathcal{Y}}  \pbr{\frac{32\sqrt{2dh_R} }{7np}\opnorm{\sum_{k\neq j} A_{jk} \Delta^{(j)}_k}  \geq \rho \frac{h_R}{\sqrt{n}} \Bigg| \Delta^{(j)}}\\
&\leq  \sup_{\Delta^{(j)}\in\mathcal{Y}} \pbr{\opnorm{\sum_{k\neq j} A_{jk} \Delta^{(j)}_k}  \geq \rho \frac{7\sqrt{n}p}{32\sqrt{d}} \Bigg| \Delta^{(j)}}, \numberthis \label{eqn:28}
\end{align*}
where in the last inequality, we use the fact that $h_R\geq 2$. Assume $\rho$ satisfies
\begin{align}
\rho\frac{7\sqrt{n}p}{64\sqrt{d}} \geq 9C_0(1+\sigma\sqrt{d}) \sqrt{dp} \label{eqn:rho_2}\,.
\end{align}
Then by Lemma \ref{lem:Fj1}, we have
\begin{align*}
&\E\indic{4\sqrt{2dh_R} \opnorm{F_{j1}} \geq \rho \frac{h_R}{\sqrt{n}}}\indic{\mathf} \\
&\leq \sup_{\Delta^{(j)}\in\mathcal{Y}} \pbr{\opnorm{\sum_{k\neq j} A_{jk} \Delta^{(j)}_k}  \geq p\sqrt{n}\opnorm{\Delta^{(j)}}+ \rho \frac{7\sqrt{n}p}{64\sqrt{d}} \Bigg| \Delta^{(j)}}\\
&\leq 2d\sup_{\Delta^{(j)}\in\mathcal{Y}}  \ebr{-\frac{\br{ \rho \frac{7\sqrt{n}p}{64\sqrt{d}} }^2/2}{p\sqrt{nd} \opnorm{\Delta^{(j)}}\opnorminf{\Delta^{(j)}}  +  \rho \frac{7\sqrt{n}p}{32\sqrt{d}} \opnorminf{\Delta^{(j)}} /3}}\\
&\leq  2d\sup_{\Delta^{(j)}\in\mathcal{Y}}  \ebr{ - \frac{3}{8} \frac{ \rho \frac{7\sqrt{n}p}{64\sqrt{d}}}{ \opnorminf{\Delta^{(j)}}}}\\
&\leq 2d \ebr{-\frac{3}{512} \frac{\rho\sqrt{n}p}{\sqrt{d}\br{\frac{ C''}{\sqrt{n}}\br{\sqrt{\frac{\log n }{np}} + \frac{\sigma \sqrt{d}}{\sqrt{np}} + \frac{\sigma \sqrt{\log n}}{\sqrt{np}}} + \frac{1}{\sqrt{n}}}}},\numberthis \label{eqn:31}
\end{align*}
where in the first and third inequality, we have used (\ref{eqn:rho_2}) and the upper bound on $\opnorm{\Delta^{(j)}}$ per $\mathcal{Y}$, and in the last inequality, we have used the upper bound on $\opnorminf{\Delta^{(j)}}$ per $\mathcal{Y}$.

\textbf{For the second term in  (\ref{eqn:hamming-decomp-ineq})}, by (\ref{eqn:3}), we have
\begin{align*}
&\E\indic{4\sqrt{2dh_R} \opnorm{F_{j2}}  \geq \rho \frac{h_R}{\sqrt{n}}}\indic{\mathf} \\
&\leq  \E\indic{ 4 \sigma \sqrt{2dh_R} \opnorm{\sum_{k\neq j} A_{jk} W_{jk}\Delta^{(j)}_k} \opnorm{\Lambda^{-1}}  \geq \rho \frac{h_R}{\sqrt{n}}}\indic{\mathf}\\
&\leq  \E\indic{\frac{32 \sigma \sqrt{2dh_R} }{7np}\opnorm{\sum_{k\neq j} A_{jk} W_{jk}\Delta^{(j)}_k}  \geq \rho \frac{h_R}{\sqrt{n}}}\indic{\mathf}\\
&\leq  \E\indic{\frac{32 \sigma \sqrt{2dh_R} }{7np}\opnorm{\sum_{k\neq j} A_{jk} W_{jk}\Delta^{(j)}_k}  \geq \rho \frac{h_R}{\sqrt{n}}}\indic{\mathf_j}.
\end{align*}
By the independence among $\{A_{jk}\}_{k\neq j}$, $\{W_{jk}\}_{k\neq j}$, and $\Delta^{(j)}$, following the same argument as used in (\ref{eqn:28}),  we have
\begin{align*}
& \E\indic{4\sqrt{2dh_R} \opnorm{F_{j2}}  \geq \rho \frac{h_R}{\sqrt{n}}}\indic{\mathf} \\
& = \E \br{\pbr{\frac{32 \sigma \sqrt{2dh_R} }{7np}\opnorm{\sum_{k\neq j} A_{jk} W_{jk}\Delta^{(j)}_k}  \geq \rho \frac{h_R}{\sqrt{n}},\mathf_j \Bigg|  \{A_{jk}\}_{k\neq j}, \Delta^{(j)}}} \\
& =  \E \br{ \pbr{\frac{32 \sigma \sqrt{2dh_R} }{7np}\opnorm{\sum_{k\neq j} A_{jk} W_{jk}\Delta^{(j)}_k}  \geq \rho \frac{h_R}{\sqrt{n}} \Bigg|  \{A_{jk}\}_{k\neq j}, \Delta^{(j)}} \indic{\Delta^{(j)}\in\mathf_j}}\\
&\leq \E\br{ \pbr{\frac{32 \sigma \sqrt{2dh_R} }{7np}\opnorm{\sum_{k\neq j} A_{jk} W_{jk}\Delta^{(j)}_k}  \geq \rho \frac{h_R}{\sqrt{n}} \Bigg| \{A_{jk}\}_{k\neq j}, \Delta^{(j)}} \indic{\Delta^{(j)}\in\mathf_j} }\\
&\leq  \E\br{ \pbr{ \opnorm{\sum_{k\neq j} A_{jk} W_{jk}\Delta^{(j)}_k}  \geq \rho \frac{7 \sqrt{n}p}{32\sigma\sqrt{d}} \Bigg| \{A_{jk}\}_{k\neq j},\Delta^{(j)} } \indic{\Delta^{(j)}\in\mathf_j}  }\\
&\leq  \E \Bigg( \br{ 2\ebr{- \frac{c\br{\rho \frac{7 \sqrt{n}p}{32\sigma\sqrt{d}} }^2}{\opnorm{ \sum_{k\neq j}A_{jk}(\Delta^{(j)}_k)^\top \Delta^{(j)}_k}}}   +  \indic{\frac{\rho \frac{7 \sqrt{n}p}{32\sigma\sqrt{d}}}{\opnorm{ \sum_{k\neq j}A_{jk}(\Delta^{(j)}_k)^\top \Delta^{(j)}_k}^\frac{1}{2}} <4\sqrt{d}} } \\
&\quad  \cdot\indic{\Delta^{(j)}\in\mathf_j} \Bigg),
\end{align*}
for some constant $c>0$, 
where the second to last inequality is by the fact that ${h_R}\geq 2$ and the last inequality is by Lemma \ref{lem:Fj2}. Since 
\begin{align*}
\ebr{- \frac{c\br{\rho \frac{7 \sqrt{n}p}{32\sigma\sqrt{d}} }^2}{\opnorm{ \sum_{k\neq j}A_{jk}(\Delta^{(j)}_k)^\top \Delta^{(j)}_k}}} &\leq \ebr{-\frac{np}{\sigma^2}} + \indic{\frac{c\br{\rho \frac{7 \sqrt{n}p}{32\sigma\sqrt{d}} }^2}{\opnorm{ \sum_{k\neq j}A_{jk}(\Delta^{(j)}_k)^\top \Delta^{(j)}_k}} \leq \frac{np}{\sigma^2}},
\end{align*}
we have
\begin{align*}
& \E\indic{2\sqrt{2dh_R} \opnorm{F_{j2}}  \geq \rho \frac{h_R}{\sqrt{n}}}\indic{\mathf} \\
&\leq  \E \left(  \br{2\ebr{-\frac{np}{\sigma^2}} +  3\indic{\frac{\br{\rho \frac{7 \sqrt{n}p}{32\sigma\sqrt{d}} }^2}{\opnorm{ \sum_{k\neq j}A_{jk}(\Delta^{(j)}_k)^\top \Delta^{(j)}_k}} \leq \frac{np}{c\sigma^2} \wedge 16d} }\indic{\Delta^{(j)}\in\mathf_j} \right)\\
&\leq   2\ebr{-\frac{np}{\sigma^2}} +  3 \E  \indic{\frac{\br{\rho \frac{7 \sqrt{n}p}{32\sigma\sqrt{d}} }^2}{\opnorm{ \sum_{k\neq j}A_{jk}(\Delta^{(j)}_k)^\top \Delta^{(j)}_k}} \leq \frac{np}{c\sigma^2} }\indic{\Delta^{(j)}\in\mathf_j}  ,
\end{align*}
where the last inequality holds as long as $\frac{np}{\sigma^2d}\geq 16c$. Then characterizing the above display is about controlling the tail probabilities of $\opnormt{ \sum_{k\neq j}A_{jk}(\Delta^{(j)}_k)^\top \Delta^{(j)}_k}$. Similar to the establishment of (\ref{eqn:28}), we have
\begin{align*}
& \E\indic{4\sqrt{2dh_R} \opnorm{F_{j2}}  \geq \rho \frac{h_R}{\sqrt{n}}}\indic{\mathf} \\
&\leq   2\ebr{-\frac{np}{\sigma^2}} +  3 \sup_{\Delta^{(j)}\in \mathcal{Y}} \pbr{\opnorm{ \sum_{k\neq j}A_{jk}(\Delta^{(j)}_k)^\top \Delta^{(j)}_k} \geq  \frac{7^2\rho^2 cp}{32^2d}\Bigg|  \Delta^{(j)}}.
\end{align*}
Assume $\rho$ satisfies
\begin{align}
& \frac{1}{2}\frac{7^2\rho^2 cp}{32^2d} \geq p\br{ \frac{9C_0(1+\sigma\sqrt{d})}{\sqrt{np}}}^2. \label{eqn:rho_3} 
\end{align}
By Lemma \ref{lem:Fj2_denominator}, we have
\begin{align*}
& \E\indic{4\sqrt{2dh_R} \opnorm{F_{j2}}  \geq \rho \frac{h_R}{\sqrt{n}}}\indic{\mathf} \\
&\leq   2\ebr{-\frac{np}{\sigma^2}} + 6d   \sup_{\Delta^{(j)}\in \mathcal{Y}} \ebr{\opnorm{ \sum_{k\neq j}A_{jk}(\Delta^{(j)}_k)^\top \Delta^{(j)}_k} \geq  p \opnorm{\Delta^{(j)}}^2 + \frac{1}{2}\frac{7^2\rho^2 cp}{32^2d} }\\
&\leq   2\ebr{-\frac{np}{\sigma^2}} + 6d   \sup_{\Delta^{(j)}\in \mathcal{Y}}  \ebr{ - \frac{\br{ \frac{1}{2}\frac{7^2\rho^2 cp}{32^2d} }^2/2}{ p \opnorminf{\Delta^{(j)}}^2\opnorm{\Delta^{(j)}}^2 + \opnorminf{\Delta^{(j)}}^2  \frac{1}{2}\frac{7^2\rho^2 cp}{16^2d} /3}}\\
&\leq   2\ebr{-\frac{np}{\sigma^2}} + 6d   \sup_{\Delta^{(j)}\in \mathcal{Y}}  \ebr{-\frac{3}{8} \frac{ \frac{1}{2}\frac{7^2\rho^2 cp}{32^2d}}{\opnorminf{\Delta^{(j)}}^2} } \\
&\leq   2\ebr{-\frac{np}{\sigma^2}} + 6d \ebr{-\frac{3}{8} \frac{ \frac{1}{2}\frac{7^2\rho^2 cp}{32^2d}}{7^2\br{\frac{ C''}{\sqrt{n}}\br{\sqrt{\frac{\log n }{np}} + \frac{\sigma \sqrt{d}}{\sqrt{np}} + \frac{\sigma \sqrt{\log n}}{\sqrt{np}}} + \frac{1}{\sqrt{n}}}^2} }\\
&\leq   2\ebr{-\frac{np}{\sigma^2}} + 6d \ebr{ -\frac{3c }{16\times 32^2} \br{ \frac{\rho \sqrt{p}}{\sqrt{d} \br{\frac{ C''}{\sqrt{n}}\br{\sqrt{\frac{\log n }{np}} + \frac{\sigma \sqrt{d}}{\sqrt{np}} + \frac{\sigma \sqrt{\log n}}{\sqrt{np}}} + \frac{1}{\sqrt{n}}}}}^2}.\numberthis \label{eqn:32}
\end{align*}

\textbf{For the third term in (\ref{eqn:hamming-decomp-ineq})}, using (\ref{eqn:3}) and the fact that $h_R\geq 2$, we have
\begin{align*}
& \E \indic{ 46C_0\br{\frac{1+\sigma \sqrt{d}}{\sqrt{np}}} \sqrt{2dh_R} \opnorm{G_{j2}} \geq  \rho \frac{h_R}{\sqrt{n}}}\indic{\mathf}  \\
&\leq \E \indic{ 46C_0\br{\frac{1+\sigma \sqrt{d}}{\sqrt{np}}} \sqrt{2dh_R} \frac{\sigma}{\sqrt{n}} \opnorm{\Lambda^{-1}}\opnorm{\sum_{k\neq j} A_{jk}W_{jk}} \geq  \rho \frac{h_R}{\sqrt{n}}}\\
 & \leq \E \indic{\opnorm{\sum_{k\neq j} A_{jk}W_{jk}} \geq  \frac{7}{8\times 46C_0\sqrt d\br{\frac{1+\sigma \sqrt{d}}{\sqrt{np}}}  } \frac{\rho np}{\sigma} }\\
 & = \E \br{\pbr{\opnorm{\sum_{k\neq j} A_{jk}W_{jk}} \geq  \frac{7}{8\times 46C_0\sqrt d\br{\frac{1+\sigma \sqrt{d}}{\sqrt{np}}}  }\frac{\rho np}{\sigma} \Bigg| \{A_{jk}\}_{k\neq j}}   }.
\end{align*}
By  the independence between $\{A_{jk}\}_{k\neq j}$ and $\{W_{jk}\}_{k\neq j}$, we have $\sum_{k\neq j} A_{jk}W_{jk} | \{A_{jk}\}_{k\neq j} \dequal \sqrt{\{A_{jk}\}_{k\neq j}} \xi$ for some $\xi \sim \MN(\zero,I_d,I_d)$ that is independent of $\{A_{jk}\}_{k\neq j}$
Using Lemma \ref{lem:vershynin_new}, we have 
\begin{align*}
& \E \indic{ 46C_0\br{\frac{1+\sigma \sqrt{d}}{\sqrt{np}}} \sqrt{2dh_R} \opnorm{G_{j2}} \geq  \rho \frac{h_R}{\sqrt{n}}}\indic{\mathf}  \\
&\leq \E \Bigg(  2\ebr{-c \br{\frac{1}{2}  \frac{7}{8\times 46C_0\sqrt d\br{\frac{1+\sigma \sqrt{d}}{\sqrt{np}}}}\frac{\rho np}{\sigma } }^2 \frac{1}{{\sum_{k\neq j}A_{jk}}} }  \\
&\quad +   \indic{\frac{1}{2}  \frac{7}{8\times 46C_0\sqrt d\br{\frac{1+\sigma \sqrt{d}}{\sqrt{np}}}  }\frac{\rho np}{\sigma} \frac{1}{\sqrt{\sum_{k\neq j}A_{jk}}} < 2\sqrt{d} }\Bigg),
\end{align*}
for some constant $c'>0$.
Note that by Bernstein's inequality, there exists some constant $C_3>0$ such that
\begin{align}\label{eqn:29}
\pbr{\sum_{k\neq j} A_{jk} \geq np + C_3\sqrt{np\log n}}\leq n^{-10}.
\end{align}
We  assume $\frac{np}{\log n}\geq C_3^2$ such that $2np\geq np + C_3\sqrt{np\log n}$ and consequently $\pbr{\sum_{k\neq j} A_{jk} \geq 2np}\leq n^{-10}$.
Assume $\rho$ satisfies
\begin{align}
\frac{1}{ \sqrt{2np}}\br{\frac{1}{2} \frac{7}{8\times 46C_0\sqrt d\br{\frac{1+\sigma \sqrt{d}}{\sqrt{np}}}  }\frac{\rho np}{\sigma}} \geq  (c')^{-\frac{1}{2}}\vee\br{2\sqrt{d}}.\label{eqn:rho_4} 
\end{align}
We have 
\begin{align*}
 \E \indic{ 46C_0\br{\frac{1+\sigma \sqrt{d}}{\sqrt{np}}} \sqrt{2dh_R} \opnorm{G_{j2}} \geq  \rho \frac{h_R}{\sqrt{n}}}\indic{\mathf}  &\leq 2\ebr{-\frac{np}{\sigma^2}} + 2\pbr{\sum_{k\neq j} A_{jk} \geq 2np}\\
&\leq   2\ebr{-\frac{np}{\sigma^2}} + 2n^{-10}.\numberthis \label{eqn:33}
\end{align*}

\textbf{For the fourth and last term in (\ref{eqn:hamming-decomp-ineq})}, by properties of the Gaussian distribution, one can verify that $\iprod{\sum_{k\neq j} A_{jk} W_{jk} }{R-I_d} \Big | \{A_{jk}\}_{k\neq j}  \sim \mathn(0,4\sum_{k\neq j} A_{jk}) $. Assume $\rho$ satisfies
\begin{align}\label{eqn:rho_5} 
3\rho + C_2\sqrt d\br{\frac{\sqrt{\log n}+ \sigma\sqrt{d}}{\sqrt{np}}} \leq  \frac{1}{2}.
\end{align}
By the fact $h_R\geq 2$, we then have
\begin{align*}
& \E \indic{ \frac{\sigma}{\sqrt{n}}\frac{1}{np}  \iprod{\sum_{k\neq j} A_{jk} W_{jk} }{R-I_d} \geq \br{1- 3\rho - C_2\sqrt d\br{\frac{\sqrt{\log n}+ \sigma\sqrt{d}}{\sqrt{np}}}}     \frac{h_R}{\sqrt{n}} }\indic{\mathf} \\
&\leq  \E \indic{  \iprod{\sum_{k\neq j} A_{jk} W_{jk} }{R-I_d} \geq \br{1- 3\rho - C_2\sqrt d\br{\frac{\sqrt{\log n}+ \sigma\sqrt{d}}{\sqrt{np}}}}    \frac{2np}{\sigma} }\\
& =\E \br{\pbr{  \iprod{\sum_{k\neq j} A_{jk} W_{jk} }{R-I_d} \geq \br{1- 3\rho - C_2\sqrt d\br{\frac{\sqrt{\log n}+ \sigma\sqrt{d}}{\sqrt{np}}}}    \frac{2np}{\sigma} \Bigg |\{A_{jk}\}_{k\neq j}}}\\
&\leq  \E \ebr{- \br{1- 3\rho - C_2\sqrt d\br{\frac{\sqrt{\log n}+ \sigma\sqrt{d}}{\sqrt{np}}}}_+^2 \frac{(np)^2}{2\sigma^2 \sum_{k\neq j} A_{jk}}}\\
&\leq  \ebr{- \br{1- 3\rho - C_2\sqrt d\br{\frac{\sqrt{\log n}+ \sigma\sqrt{d}}{\sqrt{np}}}} _+^2 \frac{(np)^2}{2\sigma^2 \br{np+C_3\sqrt{np\log n}  }}} \\
&\quad + \pbr{\sum_{k\neq j}A_{jk}\geq np+C_0\sqrt{np\log n} }\\
&\leq  \ebr{- \br{1- 3\rho - C_2\sqrt d\br{\frac{\sqrt{\log n}+ \sigma\sqrt{d}}{\sqrt{np}}}} _+^2 \br{1+C_3\sqrt{\frac{\log n}{np}}}^{-1} \frac{np}{2\sigma^2}} + n^{-10}, \numberthis \label{eqn:34}
\end{align*}
where the last inequality is by (\ref{eqn:29}).

\subsubsection{Step 4. Selecting $\rho$ and Putting Things Together} \label{subsubsec:5}

Plugging (\ref{eqn:31}), (\ref{eqn:32}), (\ref{eqn:33}),  and (\ref{eqn:34}) into (\ref{eqn:hamming-decomp-ineq}), we have
\begin{align*}
 &\E\indic{ \fnorm{U_j \hatH^\top - U_j^*\hat P^\top} \geq \fnorm{U_j \hatH^\top - \frac{R \hat P^\top}{\sqrt n} }   }  \indic{\mathf} \\
 &\leq  2d \ebr{-\frac{3}{512} \frac{\rho\sqrt{n}p}{\sqrt{d}\br{\frac{ C''}{\sqrt{n}}\br{\sqrt{\frac{\log n }{np}} + \frac{\sigma \sqrt{d}}{\sqrt{np}} + \frac{\sigma \sqrt{\log n}}{\sqrt{np}}} + \frac{1}{\sqrt{n}}}}} \\
 &\quad +  2\ebr{-\frac{np}{\sigma^2}} + 6d \ebr{ -\frac{3c }{16\times 32^2} \br{ \frac{\rho \sqrt{p}}{\sqrt{d} \br{\frac{ C''}{\sqrt{n}}\br{\sqrt{\frac{\log n }{np}} + \frac{\sigma \sqrt{d}}{\sqrt{np}} + \frac{\sigma \sqrt{\log n}}{\sqrt{np}}} + \frac{1}{\sqrt{n}}}}}^2}\\
 &\quad+ 2\ebr{-\frac{np}{\sigma^2}} + 2n^{-10} \\
 &\quad +  \ebr{- \br{1- 3\rho - C_2\sqrt d\br{\frac{\sqrt{\log n}+ \sigma\sqrt{d}}{\sqrt{np}}}}_+^2 \br{1+C_3\sqrt{\frac{\log n}{np}}}^{-1} \frac{np}{2\sigma^2}} + n^{-10}\\
 &\leq  \ebr{- \br{1- 3\rho - C_2\sqrt d\br{\frac{\sqrt{\log n}+ \sigma\sqrt{d}}{\sqrt{np}}}}_+^2 \br{1+C_3\sqrt{\frac{\log n}{np}}}^{-1} \frac{np}{2\sigma^2}} + 3n^{-10} \\
 &\quad + 4\ebr{-\frac{np}{\sigma^2}} + 2d \ebr{-\frac{3}{512} \frac{\rho\sqrt{n}p}{\sqrt{d}\br{\frac{ C''}{\sqrt{n}}\br{\sqrt{\frac{\log n }{np}} + \frac{\sigma \sqrt{d}}{\sqrt{np}} + \frac{\sigma \sqrt{\log n}}{\sqrt{np}}} + \frac{1}{\sqrt{n}}}}} \\
 &\quad + 6d \ebr{ -\frac{3c }{16\times 32^2} \br{ \frac{\rho \sqrt{p}}{\sqrt{d} \br{\frac{ C''}{\sqrt{n}}\br{\sqrt{\frac{\log n }{np}} + \frac{\sigma \sqrt{d}}{\sqrt{np}} + \frac{\sigma \sqrt{\log n}}{\sqrt{np}}} + \frac{1}{\sqrt{n}}}}}^2}. \numberthis\label{eqn:new_6}
\end{align*}
We will now pick $\rho$ such that we can substantially simplify the exponential terms. Specifically, we let $\rho = \br{\frac{\sigma^2 d^3}{np}}^\frac{1}{4} + \br{\frac{d^2\log n}{np}}^\frac{1}{4}$. Suppose that $\frac{np}{d^2\log n}, \frac{np}{\sigma^2d^3}\geq C_3$ for some large enough constant $C_3>0$, we can make sure that conditions (\ref{eqn:rho_2}), (\ref{eqn:rho_3}), (\ref{eqn:rho_4}), and (\ref{eqn:rho_5}) are all satisfied. Next, we are going to show that with this choice of $\rho$, the last two terms in (\ref{eqn:new_6}) are small. We can show the exponent in the second to last term is bounded below by $\log n$, up to a constant factor, as follows.
\begin{align*}
&\frac{\log n}{\log n} \cdot \frac{3}{512} \frac{\rho\sqrt{n}p}{\sqrt{d}\br{\frac{ C''}{\sqrt{n}}\br{\sqrt{\frac{\log n }{np}} + \frac{\sigma \sqrt{d}}{\sqrt{np}} + \frac{\sigma \sqrt{\log n}}{\sqrt{np}}} + \frac{1}{\sqrt{n}}}} \\
&= \log n \cdot \frac{3}{512}\frac{\rho}{ \frac{\sqrt d \log n}{np} \br{C''\sqrt{\frac{\log n}{np}} + \frac{C'' \sigma \sqrt d}{\sqrt{np}} + \frac{C''\sigma\sqrt{\log n}}{\sqrt{np}} + 1   }  }\\
&\geq  \log n \cdot \frac{3}{512}\frac{\rho}{ \frac{\sqrt d (\log n)^{3/2}}{np} \br{C''\sqrt{\frac{\log n}{np}} + \frac{C'' \sigma \sqrt d}{\sqrt{np}} + \frac{C''\sigma}{\sqrt{np}} + 1   }  }\\
&\geq  \log n \cdot \frac{3}{512}\frac{\br{\frac{d^2\log n}{np}}^\frac{1}{4}}{ \frac{\sqrt d (\log n)^{3/2}}{np} \br{C''\sqrt{\frac{\log n}{np}} + \frac{C'' \sigma \sqrt d}{\sqrt{np}} + \frac{C''\sigma}{\sqrt{np}} + 1   }  }\\
&=  \log n \cdot \frac{3}{512}\frac{1}{ \frac{(\log n)^{5/4}}{(np)^{3/4}} \br{C''\sqrt{\frac{\log n}{np}} + \frac{C'' \sigma \sqrt d}{\sqrt{np}} + \frac{C''\sigma}{\sqrt{np}} + 1   }  }\\
&\geq 10\log n.
\end{align*}
In the first line, we simply multiply and divide the same quantity by $\log n$. The second and third lines are by expanding and re-arranging the terms. In the second inequality, we use $\rho \geq \br{\frac{d^2\log n}{np}}^{\frac{1}{4}}$. The last inequality holds so long as $\frac{np}{d\log n}, \frac{np}{\log^{5/3}n}, \frac{np}{\sigma^2 d} \geq C_3'$ for some sufficiently large constant $C_3'$.

Similarly, we can also simplify the term inside the square root of the last term in (\ref{eqn:new_6}) as
\begin{align*}
&\frac{\sqrt{\log n}}{\sqrt{\log n}} \cdot \frac{\sqrt{3c}}{4\times 32} \cdot  \frac{\rho \sqrt{p}}{\sqrt{d} \br{\frac{ C''}{\sqrt{n}}\br{\sqrt{\frac{\log n }{np}} + \frac{\sigma \sqrt{d}}{\sqrt{np}} + \frac{\sigma \sqrt{\log n}}{\sqrt{np}}} + \frac{1}{\sqrt{n}}}}\\
&=\sqrt{\log n} \cdot \frac{\sqrt{3c}}{4\times 32} \cdot  \frac{\rho }{\frac{\sqrt{d\log n}}{\sqrt{np}} \br{C''\sqrt{\frac{\log n }{np}} +C'' \frac{\sigma \sqrt{d}}{\sqrt{np}} +C'' \frac{\sigma \sqrt{\log n}}{\sqrt{np}} + 1}}\\
&\geq \sqrt{\log n} \cdot \frac{\sqrt{3c}}{4\times 32} \cdot  \frac{\rho }{\frac{\sqrt{d} \log n}{\sqrt{np}} \br{C''\sqrt{\frac{\log n }{np}} +C'' \frac{\sigma \sqrt{d}}{\sqrt{np}} +C'' \frac{\sigma }{\sqrt{np}} + 1}}\\
&\geq \sqrt{\log n} \cdot \frac{\sqrt{3c}}{4\times 32} \cdot  \frac{ \br{\frac{d^2\log n}{np}}^{\frac{1}{4}}}{\frac{\sqrt{d} \log n}{\sqrt{np}} \br{C''\sqrt{\frac{\log n }{np}} +C'' \frac{\sigma \sqrt{d}}{\sqrt{np}} +C'' \frac{\sigma }{\sqrt{np}} + 1}}\\
&\geq \sqrt{\log n} \cdot \frac{\sqrt{3c}}{4\times 32} \cdot  \frac{1}{\frac{(\log n)^{3/4}}{(np)^{1/4}} \br{C''\sqrt{\frac{\log n }{np}} +C'' \frac{\sigma \sqrt{d}}{\sqrt{np}} +C'' \frac{\sigma }{\sqrt{np}} + 1}}\\
&\geq \sqrt{10}\sqrt{\log n}.
\end{align*}
In the first line, we multiply and divide the same quantity by $\sqrt{\log n}$. The last inequality holds so long as $\frac{np}{d\log n},\frac{np}{\log^3 n}, \frac{np}{\sigma^2 d} \geq C_3''$ for some sufficiently large constant $C_3''$.

Putting things together into (\ref{eqn:new_6}), we have
\begin{align*}
&\E\indic{ \fnorm{U_j \hatH^\top - U_j^*\hat P^\top} \geq \fnorm{U_j \hatH^\top - \frac{R \hat P^\top}{\sqrt n} }   }  \indic{\mathf} \\
&\leq   \exp\Bigg(- \br{1- 3\br{\br{\frac{\sigma^2d^3}{np}}^\frac{1}{4} + \br{\frac{d^2\log n}{np}}^\frac{1}{4}} - C_2\sqrt d\br{\frac{\sqrt{\log n}+ \sigma\sqrt{d}}{\sqrt{np}}}}_+^2 \\
&\quad \cdot \br{1+C_3\sqrt{\frac{\log n}{np}}}^{-1} \frac{np}{2\sigma^2}\Bigg)\\
&\quad  + 4\ebr{-\frac{np}{\sigma^2}} + 3n^{-10} + \underbrace{2d\ebr{-10\log n} + 6d\ebr{-10\log n}}_{\leq 9dn^{-10}}\,,
\end{align*}
By combining the two exponential terms together and the three polynomial terms together, we see that there exists a constant $C_4>0$ such that
\begin{align*}
&\E\indic{ \fnorm{U_j \hatH^\top - U_j^*\hat P^\top} \geq \fnorm{U_j \hatH^\top - \frac{R \hat P^\top}{\sqrt n} }   }  \indic{\mathf} \\
&\leq \ebr{-\br{1-C_4\br{\br{\frac{\sigma^2d^3}{np}}^\frac{1}{4} + \br{\frac{d^2\log n}{np}}^\frac{1}{4} }}_+\frac{np}{2\sigma^2}} + 12dn^{-10}. \numberthis\label{eqn:tail-bound}
\end{align*}

Combining (\ref{eqn:36}), (\ref{eqn:35}), (\ref{eqn:26}) and (\ref{eqn:tail-bound}) we have
\begin{align*}
&\E \ell(\hat Z,Z^*) \\
&\leq \frac{1}{n}\sum_{j\in[n]}\sum_{R\in \Pi_d:R\neq I_d} \br{\ebr{-\br{1-C_4\br{\br{\frac{\sigma^2d^3}{np}}^\frac{1}{4} + \br{\frac{d^2\log n}{ np}}^\frac{1}{4} }}_+\frac{np}{2\sigma^2}} + 12dn^{-10}} + \pbr{\mathf^\mathrm{c}}\\
&\leq d^d \br{\ebr{-\br{1-C_4\br{\br{\frac{\sigma^2d^3}{np}}^\frac{1}{4} + \br{\frac{d^2\log n}{np}}^\frac{1}{4} }}_+\frac{np}{2\sigma^2}} + 12dn^{-10}} + 7dn^{-9}\\
&\leq \ebr{-\br{1-C_5\br{\br{\frac{\sigma^2d^3}{np}}^\frac{1}{4} + \br{\frac{d^2\log n}{np}}^\frac{1}{4} }}_+\frac{np}{2\sigma^2}}  + 8dn^{-9}\\
&\leq \ebr{-\br{1-C_5\br{\br{\frac{\sigma^2d^3}{np}}^\frac{1}{4} + \br{\frac{d^2\log n}{np}}^\frac{1}{4} }}\frac{np}{2\sigma^2}}  + 8dn^{-9},
\end{align*}
for some constant $C_5>0$, where the last inequality holds when $\frac{np}{\sigma^2d^3}, \frac{np}{d^2\log n}$ exceed a sufficiently large constant. The proof of Theorem \ref{thm:main_new} is complete.

\subsubsection{Proof of Theorem \ref{thm:main_approx}.}  \label{subsubsec:6}

The proof closely parallels that of Theorem \ref{thm:main_new}, with substitutions where $\hat{Z}$ and $M$ are replaced by $\check{Z}$ and $\check{M}$, respectively, and incorporates two minor modifications. Firstly, it is important to note that the error in the anchor affects only the second inner product in (\ref{eqn:9}). Specifically, this influence introduces an additional factor of $(1+\epsilon)$, leading to a modified formulation in  (\ref{eqn:second-inner}) for the scenario involving $\check{M}$:
\begin{align*}
&\iprod{U_j\br{\check M - \hat PH^\top}^\top \hat P}{ R-I_d} \\
&\leq  \fnorm{U_j \br{\check M -\hat P  H^\top}^\top}\fnorm{ R-I_d}\\
&\leq  \opnorm{U_j} \fnorm{\check M -\hat P  H^\top}\fnorm{ R-I_d}\\
&\leq   \frac{4}{3}\br{\opnorm{U_jH-U^*_j}+\frac{1}{\sqrt{n}}}\br{ \frac{{4C_0\sqrt{1+\epsilon}\br{\sqrt{d} + \sigma d}}}{\sqrt{np}}}\sqrt{2h_R}\\
&\leq  \frac{4}{3}\br{\opnorm{B_j}+\opnorm{F_{j1}}+\opnorm{F_{j2}}+\opnorm{G_{j1}}+ \opnorm{G_{j2}}+\frac{1}{\sqrt{n}}}\br{ \frac{{4C_0\sqrt{1+\epsilon}\br{\sqrt{d} + \sigma d}}}{\sqrt{np}}}\sqrt{2h_R},
\end{align*}
where in the second inequality, we use (\ref{eqn:check_M}) from Lemma \ref{lem:check_M} instead (\ref{eqn:11}) from Lemma \ref{lem:hat_H}.
We can then follow the same argument in the proof leading up to (\ref{eqn:25}) where we make the second modification.
\begin{align*}
&\iprod{U_j H - U_j^*}{R-I_d} + \iprod{U_j\br{\check M - \hat PH^\top}^\top \hat P}{ R-I_d} \\
&\leq  2\br{1+ \frac{16C_0\sqrt{1+\epsilon}(1+\sigma\sqrt{d})}{3\sqrt{np}} }\br{ \opnorm{F_{j1}} + \opnorm{F_{j2}} } \sqrt{2dh_R} \\
&\quad +\br{\frac{41}{33} +\frac{16\sqrt{1+\epsilon}C_0(1+\sigma\sqrt{d})}{3\sqrt{np}}  }23C_0\br{\frac{1+\sigma \sqrt{d}}{\sqrt{np}}}  \opnorm{G_{j2}} \sqrt{2dh_R} \\
&\quad +\br{\frac{41}{33} +\frac{16\sqrt{1+\epsilon}C_0(1+\sigma\sqrt{d})}{3\sqrt{np}}  }26C_0\br{\frac{\sqrt{\log n}+\sigma \sqrt{d}}{\sqrt{np}}} \frac{\sqrt{2dh_R} }{\sqrt{n}}\\
&\quad +  \frac{\sigma}{\sqrt{n}}\frac{1}{np}  \iprod{\sum_{k\neq j} A_{jk} W_{jk} }{R-I_d}\\
&\leq 4\br{ \opnorm{F_{j1}} + \opnorm{F_{j2}} } \sqrt{2dh_R}  + 46C_0\br{\frac{1+\sigma \sqrt{d}}{\sqrt{np}}}  \opnorm{G_{j2}} \sqrt{2dh_R} \\
&\quad + 52C_0\br{\frac{\sqrt{\log n}+\sigma \sqrt{d}}{\sqrt{np}}} \frac{\sqrt{2dh_R} }{\sqrt{n}} +  \frac{\sigma}{\sqrt{n}}\frac{1}{np}  \iprod{\sum_{k\neq j} A_{jk} W_{jk} }{R-I_d},
\end{align*}
where the last inequality holds under the  assumption that $\frac{\sqrt{np}}{1+\sigma \sqrt d} > \frac{176}{25}C_0\sqrt{1+\epsilon}$. In the way, we obtain the exact same upper bound as (\ref{eqn:25}). There is no change to the rest of the proof.

\subsection{Helper Lemmas and Proofs}
\begin{lemma}\label{lem:population}
We have $\lambda_i(\E A\otimes I_d)=(n-1)p$ for all $i\leq d$ and $\lambda_i(\E A\otimes I_d)=-p$ for all $d+1\leq i\leq nd$. 

\end{lemma}
\begin{proof} Recall that $\E A = pJ_n - pI_n$. The eigenvalues of the matrix $\E A$ are characterized as follows.
$$ \lambda_1(\E A) = (n-1)p \,,$$
$$ \lambda_2(\E A) = \ldots = \lambda_n(\E A) = -p\,. $$
Therefore the eigenvalues of $\E A \otimes I_d$ are as follows.
$$ \lambda_1(\E A\otimes I_d) = \ldots = \lambda_d(\E A \otimes I_d) = (n-1)p \,,$$
$$\lambda_{d+1}(\E A\otimes I_d) = \ldots = \lambda_{dn}(\E A \otimes I_d) = -p\,. $$
\end{proof}

\begin{delayedproof}{lem:event-E-prime_new}
The probability of $\calE$ is from Lemma 9 of \cite{zhang2022exact} (for $ \opnorm{(A \otimes J_d)\circ W} $), Theorem 5.2 of \cite{lei2015consistency} (for $\opnorm{A - \E A} $), and Bernstein's inequality (for $\max_{j\in[n]}|\sum_{k\neq j}A_{jk}-np|$). Note that
\begin{align*}
\opnorm{X - (\E A \otimes I_d)}  & = \opnorm{ (A \otimes I_d) + \sigma (A \otimes J_d) \circ W - (\E A \otimes I_d)} \\
&\leq \opnorm{(A-\E A) \otimes I_d} +\sigma\opnorm{(A \otimes J_d) \circ W }\\
& = \opnorm{A -\E A } +\sigma\opnorm{(A \otimes J_d) \circ W }.
\end{align*}
By Weyl's inequality, we have $\max_{i\in[nd]}\abs{\lambda_i(X) - \lambda_i(\E A\otimes I_d)}\leq \opnorm{X - (\E A \otimes J_d)} $. Together with Lemma \ref{lem:population}, (\ref{eqn:2}) holds under the event $\calE$. For (\ref{eqn:7}), we have
\begin{align*}
\max_j\opnorm{X_j} &= \max_j \opnorm{  (\E A_j \otimes I_d)+ \br{X_j - (\E A_j \otimes I_d)}} \\
&\leq  \max_j \opnorm{  (\E A_j \otimes I_d)} + \opnorm{X - (\E A \otimes I_d)}\\
&\leq p\sqrt{n} + C_0\br{1+\sigma\sqrt{d}}\sqrt{np}.
\end{align*}

If $\frac{np}{(1+\sigma \sqrt{d})^2 }\geq 64C_0^2$ is further assumed, we have $C_0({1+\sigma\sqrt{d}})\sqrt{np} \leq np/8$. Then (\ref{eqn:2}), together with the fact $\opnorm{\Lambda^{-1}} = 1/\lambda_d(X)$ and $\opnorm{\Lambda} = \lambda_1(X)$, leads to (\ref{eqn:8})-(\ref{eqn:4}).

Since $\lambda_d(\E A \otimes I_d) - \lambda_{d+1}(\E A \otimes I_d) = np$ according to Lemma \ref{lem:population}, we have $\opnorm{X - (\E A \otimes I_d)}  \leq (\lambda_d(\E A \otimes I_d) - \lambda_{d+1}(\E A \otimes I_d))/8$.  
By Lemma 2 of \cite{abbe2020entrywise}, we have
\begin{align*}
\min_{O\in \mathcal{O}_d}\opnorm{H - O} \leq  \opnorm{UU^\top - U^* U^{*\top} }\leq \frac{\opnorm{X - (\E A \otimes I_d)} }{\br{1-\frac{1}{8}}np}\leq  \frac{8C_0\br{1+\sigma\sqrt{d}}}{7\sqrt{np}},
\end{align*}
and
\begin{align*}
\opnorm{H\Lambda - \Lambda H}\leq 2\opnorm{X - (\E A \otimes I_d)} \leq 2C_0\br{1+\sigma\sqrt{d}}\sqrt{np},
\end{align*}
and $\opnorm{H^{-1}}\leq 4/3$. Note that 
$$\opnorm{UH-U^*} = \opnorm{UU^\top U^* - U^*U^{*\top}U^*} \leq  \opnorm{UU^\top  - U^*U^{*\top}}\opnorm{U^*}\leq \opnorm{UU^\top  - U^*U^{*\top}}.$$
Hence, (\ref{eqn:5})-(\ref{eqn:6}) hold. Consider any $j\in[n]$. By  (\ref{eqn:6}), we have 
\begin{align*}
\opnorm{U_j} &\leq \opnorm{U_j H}\opnorm{H^{-1}} \leq \br{\opnorm{U_jH-U^*_j} + \opnorm{U^*_j}}\opnorm{H^{-1}} = \frac{4}{3}\br{\opnorm{U_jH-U^*_j} +\frac{1}{\sqrt{n}} }.
\end{align*}
From (\ref{eqn:7}), (\ref{eqn:4}) and that $C_0 > 7$, we have $\opnorm{X-\loo{X}} \leq 2\opnorm{X_j} \leq (\lambda_d(X)- \lambda_{d+1}(X))/2$. 
By Davis-Kahan theorem, we have
\begin{align*}
\opnorm{\loo{U}(\loo{U})^\top -UU^\top}  &\leq \frac{2\opnorm{\br{X-\loo{X}}U}}{\lambda_d(X) - \lambda_{d+1}(X)} \\
&\leq \frac{2\br{\opnorm{X_j U} + \opnorm{X_j} \opnorm{U_j}}}{\lambda_d(X) - \lambda_{d+1}(X)}\\
& =   \frac{2\br{\opnorm{U_j\Lambda} + \opnorm{X_j} \opnorm{U_j}}}{\lambda_d(X) - \lambda_{d+1}(X)}\\
& \leq   \frac{2\br{\opnorm{U_j}\opnorm{\Lambda} + \opnorm{X_j} \opnorm{U_j}}}{\lambda_d(X) - \lambda_{d+1}(X)}\\
&\leq 4\opnorm{U_j}\\
&\leq  6\br{\opnorm{U_jH -U^*_j}+\frac{1}{\sqrt{n}}},
\end{align*}
where the second to last inequality is due to (\ref{eqn:7}), (\ref{eqn:8}), and (\ref{eqn:4}), and the last inequality is due to (\ref{eqn:10}). 

Together with (\ref{eqn:5}) and that $C_0 > 7$ we have
\begin{align*}
\opnorm{\loo{U}(\loo{U})^\top -UU^\top} &\leq   6\br{\opnorm{UH -U^*}+\frac{1}{\sqrt{n}}}\\
&\leq \frac{48C_0(1+\sigma \sqrt{d})}{7\sqrt{np}} + \frac{6}{\sqrt{np}}\\
&\leq  \frac{7C_0(1+\sigma \sqrt{d})}{\sqrt{np}}.
\end{align*}
Then,
\begin{align*}
\opnorm{\loo{U}\loo{H} - U^*} &= \opnorm{\br{\loo{U}\loo{U}^\top - UU^\top}U^* + UH- U^*}  \\
&\leq \opnorm{\loo{U}\loo{U}^\top - UU^\top} + \opnorm{UH- U^*}\\
&\leq  \frac{9C_0(1+\sigma \sqrt{d})}{\sqrt{np}}.
\end{align*}
For $\opnorm{G_{j1}}$, we have
\begin{align*}
\sqrt{n}\opnorm{G_{j1}} &= \opnorm{\sum_{k\neq j} A_{jk} \Lambda^{-1} - I_d}\\
&\leq \opnorm{np \Lambda^{-1}-I_d} + \opnorm{ \br{\sum_{k\neq j} A_{jk} - np}\Lambda^{-1}}\\
&\leq \opnorm{np -\Lambda} \opnorm{\Lambda^{-1}} + \abs{\sum_{k\neq j} A_{jk} - np}\opnorm{\Lambda^{-1}}\\
&\leq \max_{i\in[d]} \abs{np - \lambda_i(X)}\opnorm{\Lambda^{-1}} + \abs{\sum_{k\neq j} A_{jk} - np} \opnorm{\Lambda^{-1}}\\
&\leq 2C_0\br{\sqrt{\frac{\log n}{np}} + \frac{\sigma \sqrt{d}}{\sqrt{np}}},
\end{align*}
where the last inequality is due to $\calE$,  (\ref{eqn:2}), and (\ref{eqn:3}).

For (\ref{eqn:15}), consider any $j\in[n]$, we have
\begin{align*}
\opnorminf{\loo{U}\loo{H} - U^*  } &\leq \opnorminf{UH - U^*} +  \opnorminf{\br{UH - U^*} - \br{\loo{U}\loo{H} - U^* }} \\
& \leq \opnorminf{UH - U^*} +   \opnorm{\br{UH - U^*} - \br{\loo{U}\loo{H} - U^* }  } \\
&= \opnorminf{UH - U^*} + \opnorm{UH-U^{(j)}H^{(j)}} \\
&= \opnorminf{UH - U^*} + \opnorm{UU^\top U^*-U^{(j)}U^{(j)\top} U^*}\\
&\leq  \opnorminf{UH - U^*} + \opnorm{UU^\top -U^{(j)}U^{(j)\top}} \\
&\leq  \opnorminf{UH - U^*} + 6\br{\opnorm{U_jH -U^*_j}+\frac{1}{\sqrt{n}}}\\
&\leq 7 \br{ \opnorminf{UH - U^*}+\frac{1}{\sqrt{n}}}, 
\end{align*}
where the second to last inequality is due to  (\ref{eqn:14}).

\end{delayedproof}

\begin{delayedproof}{lem:Fj1}
Consider any $j\in[n]$. 
Then
\begin{align*}
\opnorm{\sum_{k\neq j} A_{jk}\Delta_k^{(j)} } &\leq  \opnorm{\sum_{k\neq j} (A_{jk}-p)\Delta_k^{(j)}} + p\opnorm{\sum_{k\neq j} \Delta_k^{(j)}}\\
&\leq \opnorm{\sum_{k\neq j} (A_{jk}-p)\Delta_k^{(j)}} + p\sqrt{n} \opnorm{\Delta^{(j)}}.
\end{align*}
Since $\{A_{jk}-p\}_{k\neq j}$ is independent of $\Delta^{(j)}$, we use the matrix Bernstein's inequality (Lemma \ref{lem:matrix-bernstein-inequality_new}) for the operator norm of $\sum_{k\neq j} (A_{jk}-p)\Delta_k^{(j)}$. For each $k\in[n]$, note that $\E \br{(A_{jk}-p)\Delta_k^{(j)}  | \Delta_k^{(j)} }=0$ and $\opnorm{(A_{jk}-p)\Delta_k^{(j)}}\leq \opnorm{\Delta_k^{(j)}}\leq \opnorminf{\Delta^{(j)}}$. For the matrix variance term, we have
\begin{align*}
&\max\left\{ \opnorm{\E \left[ \sum_{k\neq j} (A_{jk}-p)^2 {\Delta_k^{(j)}}^\top  \Delta_k^{(j)}  \Bigg| \Delta^{(j)} \right]} , \opnorm{\E \left[ \sum_{k\neq j} (A_{jk}-p)^2  \Delta_k^{(j)}    {\Delta_k^{(j)}}^\top  \Bigg| \Delta^{(j)} \right] }  \right\}\\
&= p(1-p) \cdot \max\left\{ \opnorm{\sum_{k\neq j} {\Delta_k^{(j)}} ^\top  \Delta_k^{(j)} }, \opnorm{\sum_{k\neq j} \Delta_k^{(j)}{\Delta_k^{(j)} }^\top }   \right\} \\
&= p(1-p) \cdot \max\left\{  \opnorm{\Delta^{(j)}}^2, \opnorm{\br{\loo{\Delta_1},\ldots, \loo{\Delta_n}}^\top}^2 \right\}. \numberthis \label{eqn:lem_Fj1_1}
\end{align*}

To obtain the two different tail bounds in the lemma statement, we will show two different upper bounds on (\ref{eqn:lem_Fj1_1}).
{For the first tail bound (\ref{eqn:Fj1-tail1}) of Lemma \ref{lem:Fj1},} note that $ \opnorm{\Delta^{(j)}}^2\leq \sqrt{n}\opnorm{\Delta^{(j)}}\opnorminf{\Delta^{(j)}}$. We have 
$$p(1-p) \opnorm{\loo{\Delta}}^2 \leq p\sqrt{n} \opnorminf{\loo{\Delta}} \opnorm{\loo{\Delta}} \leq p\sqrt{n}\opnorminf{\loo{\Delta}} \opnorm{\loo{\Delta}}.$$ 
On the other hand,
\begin{align*}
 p(1-p)\opnorm{\br{\loo{\Delta_1},\ldots, \loo{\Delta_n}}^\top}^2&\leq p\sqrt{n}\opnorm{\br{\loo{\Delta_1},\ldots, \loo{\Delta_n}}^\top}\opnorminf{\Delta^{(j)}}\\
 &\leq p\sqrt{n}\fnorm{\br{\loo{\Delta_1}^\top,\ldots, \loo{\Delta_n}^\top}}\opnorminf{\Delta^{(j)}}\\
 &= p\sqrt{n}\fnorm{\Delta^{(j)}}\opnorminf{\Delta^{(j)}}\\
 &\leq p\sqrt{nd}\opnorm{\Delta^{(j)}}\opnorminf{\Delta^{(j)}}.
\end{align*}
Then (\ref{eqn:lem_Fj1_1}) can be upper bounded by $p\sqrt{nd}\opnorm{\Delta^{(j)}}\opnorminf{\Delta^{(j)}}$. 
Invoking Lemma \ref{lem:matrix-bernstein-inequality_new}, we have
\begin{align*}
&\pbr{\opnorm{\sum_{k\neq j} (A_{jk} - p) \Delta_k^{(j)}}  \geq  t\Bigg| \Delta^{(j)}  }\\
&\leq 2d\ebr{- \frac{t^2/2}{p\sqrt{nd}\opnorm{\Delta^{(j)}}\opnorminf{\Delta^{(j)}} + \opnorminf{\Delta^{(j)}}t/3}}.
\end{align*}
Then, by triangle inequality,
\begin{align*}
&\pbr{\opnorm{\sum_{k\neq j} A_{jk} \Delta_k^{(j)}}  \geq  p\sqrt{n} \opnorm{\Delta^{(j)}} + t\Bigg| \Delta^{(j)} } \\
&\leq 2d\ebr{- \frac{t^2/2}{p\sqrt{nd}\opnorm{\Delta^{(j)}}\opnorminf{\Delta^{(j)}} + \opnorminf{\Delta^{(j)}}t/3}}.
\end{align*}
{For the second tail bound (\ref{eqn:Fj1-tail2}) of Lemma \ref{lem:Fj1},} one can see that 
\begin{align*}
p(1-p) \cdot \max\left\{  \opnorm{\Delta^{(j)}}^2, \opnorm{\br{\loo{\Delta_1}^\top,\ldots, \loo{\Delta_n}^\top}}^2 \right\} \leq np \opnorminf{\loo{\Delta}}^2.
\end{align*}
We then follow the same steps as shown above to obtain (\ref{eqn:Fj1-tail2}).
\end{delayedproof}

\begin{delayedproof}{lem:Fj2_denominator}
We follow the proof of Lemma \ref{lem:Fj1}. We first have
\begin{align*}
\opnorm{ \sum_{k\neq j}A_{jk}(\Delta^{(j)}_k)^\top \Delta^{(j)}_k} &\leq  \opnorm{ \sum_{k\neq j}(A_{jk}-p)(\Delta^{(j)}_k)^\top \Delta^{(j)}_k}  + p\opnorm{ \sum_{k\neq j}(\Delta^{(j)}_k)^\top \Delta^{(j)}_k} \\
&\leq  \opnorm{ \sum_{k\neq j}(A_{jk}-p)(\Delta^{(j)}_k)^\top \Delta^{(j)}_k}  + p \opnorm{\Delta^{(j)}}^2.
\end{align*}
Note that for each $k\in[n]$, $\E \br{(A_{jk}-p)(\Delta^{(j)}_k)^\top \Delta^{(j)}_k |\Delta^{(j)}} =0$ and $\opnorm{(A_{jk}-p)(\Delta^{(j)}_k)^\top \Delta^{(j)}_k}\leq \opnorm{ \Delta^{(j)}_k}^2 \leq \opnorminf{\Delta^{(j)}}^2$. In addition,
\begin{align*}
&\opnorm{\E \left[ \sum_{k\neq j} (A_{jk}-p)^2 (\Delta^{(j)}_k)^\top \Delta^{(j)}_k(\Delta^{(j)}_k)^\top \Delta^{(j)}_k  \Bigg| \Delta^{(j)} \right]} \\
& =p(1-p)\opnorm{\sum_{k\neq j}  (\Delta^{(j)}_k)^\top \Delta^{(j)}_k(\Delta^{(j)}_k)^\top \Delta^{(j)}_k }  \\
&\leq p(1-p) \opnorm{ \,\left[{\Delta^{(j)}_1}^\top \,\cdots\, {\Delta^{(j)}_n}^\top   \right]  \begin{bmatrix} {\Delta^{(j)}_1}{\Delta^{(j)}_1}^\top & &\\ & \cdots &\\  & & {\Delta^{(j)}_n}{\Delta^{(j)}_n}^\top  \end{bmatrix}  \begin{bmatrix} \Delta^{(j)}_1\\ \vdots \\ \Delta^{(j)}_n \end{bmatrix} } \\
&\leq p(1-p)  \opnorm{ \begin{bmatrix} {\Delta^{(j)}_1}{\Delta^{(j)}_1}^\top & &\\ & \cdots &\\  & & {\Delta^{(j)}_n}{\Delta^{(j)}_n}^\top  \end{bmatrix}   }  \opnorm{ \Delta^{(j)} }^2\\
&=p(1-p)  \br{\max_{k\in[n]} \opnorm{\Delta^{(j)}_k{\Delta^{(j)}_k}^\top}} \opnorm{ \Delta^{(j)} }^2\\
&\leq  p(1-p)  \opnorminf{\Delta^{(j)}}^2  \opnorm{ \Delta^{(j)} }^2 \leq p  \opnorminf{\Delta^{(j)}}^2  \opnorm{ \Delta^{(j)} }^2\,.
\end{align*}
Hence, by Lemma \ref{lem:matrix-bernstein-inequality_new},
\begin{align*}
&\pbr{\opnorm{ \sum_{k\neq j}A_{jk}(\Delta^{(j)}_k)^\top \Delta^{(j)}_k} \geq  p \opnorm{\Delta^{(j)}}^2 + t \Bigg| \Delta^{(j)} } \\
&\leq 2d\ebr{ - \frac{t^2/2}{p \opnorminf{\Delta^{(j)}}^2  \opnorm{ \Delta^{(j)} }^2 + \opnorminf{\Delta^{(j)}}^2t/3}}.
\end{align*}
\end{delayedproof}

\begin{delayedproof}{lem:Fj2}
Consider any $j\in[n]$. Note that $\{A_{jk}\}_{k\neq j}$, $\{W_{jk}\}_{k\neq j}$, and $\{\Delta^{(j)}_k\}_{k\neq j}$ are mutually independent of each other. By Lemma \ref{lem:useful-MN-identities_new}, we have
\begin{align*}
\sum_{k\neq j}A_{jk}W_{jk}\Delta^{(j)}_k \Bigg| \{A_{jk}\}_{k\neq j}, \Delta^{(j)} \dequal E \br{\sum_{k\neq j}A_{jk}(\Delta^{(j)}_k)^\top \Delta^{(j)}_k}^\frac{1}{2} \Bigg| \{A_{jk}\}_{k\neq j},\Delta^{(j)}
\end{align*}
where $E \sim \MN(\zero, I_d, I_d)$ and is independent of $\{A_{jk}\}_{k\neq j}$ and  $\{\Delta^{(j)}_k\}_{k\neq j}$. Since $\sum_{k\neq j}A_{jk}(\Delta^{(j)}_k)^\top \Delta^{(j)}_k$ is symmetric and positive semi-definite, its square root is well-defined. By Lemma \ref{lem:vershynin_new}, there exists some constant $c>0$, such that for any $t\geq 4\sqrt{d}\opnorm{ \sum_{k\neq j}A_{jk}(\Delta^{(j)}_k)^\top \Delta^{(j)}_k}^\frac{1}{2} $, we have
\begin{align*}
& \pbr{\opnorm{\sum_{k\neq j}A_{jk}W_{jk}\Delta^{(j)}_k } \geq t  \Bigg| \{A_{jk}\}_{k\neq j}, \Delta^{(j)}} \\
 &= \pbr{\opnorm{E \br{\sum_{k\neq j}A_{jk}(\Delta^{(j)}_k)^\top \Delta^{(j)}_k}^\frac{1}{2}} \geq t \Bigg| \{A_{jk}\}_{k\neq j}, \Delta^{(j)}}\\
&\leq  \pbr{\opnorm{E}\opnorm{ \br{\sum_{k\neq j}A_{jk}(\Delta^{(j)}_k)^\top \Delta^{(j)}_k}^\frac{1}{2}}   \geq t  \Bigg| \{A_{jk}\}_{k\neq j}, \Delta^{(j)}}\\
& = \pbr{\opnorm{E}\opnorm{ \sum_{k\neq j}A_{jk}(\Delta^{(j)}_k)^\top \Delta^{(j)}_k}^\frac{1}{2}   \geq t  \Bigg| \{A_{jk}\}_{k\neq j}, \Delta^{(j)}}\\
&\leq 2\ebr{- \frac{ct^2}{\opnorm{ \sum_{k\neq j}A_{jk}(\Delta^{(j)}_k)^\top \Delta^{(j)}_k}}}.
\end{align*}
\end{delayedproof}

\begin{lemma}\label{lem:matrix-bernstein-inequality_new} [Theorem 1.6 of \cite{tropp2012user}] Consider a finite sequence $\{S_k\}_{k=1}^n$ of independent random square matrices with dimension $d\times d$. Assume that each random matrix satisfies
$$\E S_k =0 \text{ and } \opnorm{S_k - \E S_k} \leq L, \,\forall k \in [n] \,.$$
Define
\begin{align*}
V &\define  \max\left\{ \opnorm{\sum_k \E\left[ S_k S_k ^\top \right]}, \opnorm{ \sum_k\E\left[ S_k^{\top}S_k \right]}  \right\} \,.
\end{align*}
Then
$$ \pbr{\opnorm{\sum_k S_k} \geq t} \leq 2d \cdot \ebr{-\frac{t^2/2}{V + Lt/3}} \,,$$
for any $t\geq 0$.
\end{lemma}

\begin{lemma}\label{lem:useful-MN-identities_new}
Consider any $W,W',\Delta \in\mathr^{d\times d}$ such that $W,W'\iid \MN(\zero, I_d, I_d)$ and $\Delta$ is a fixed matrix.
We have
$$ W\Delta \sim \MN(\zero, I_d, \Delta^\top \Delta)\, $$
and
$$ W + W' \sim \MN(\zero, I_d, \Sigma_1 + \Sigma_2) \,.$$
\end{lemma}
\begin{proof} 
The first result in the lemma is a property of the matrix normal distribution. To prove the second statement, note that $Z \sim \MN(\zero, \Sigma_0, \Sigma_1)$ is equivalent to $\vec(Z) \sim \calN(0, \Sigma_1 \otimes \Sigma_0)$. Since $ \vec(W) \sim \calN(\zero,  \Sigma_1 \otimes I_d) $ and $ \vec(W') \sim \calN(\zero,  \Sigma_2 \otimes I_d) ,$ we have
$ \vec(W) + \vec(W') \sim \calN(\zero,  (\Sigma_1 +\Sigma_2) \otimes I_d).$
\end{proof}

\bibliography{reference}

\end{document}